\documentclass[11pt]{article}
\usepackage[top=1in,bottom=1in,left=1in,right=1in]{geometry}
\usepackage{graphics,graphicx,amsmath,amssymb,amsthm,eqnarray,array,xcolor,verbatim,easytable,subcaption,hyperref,setspace,csquotes,rotating,float,nomencl,mathtools,ulem,booktabs,units,authblk,abstract,appendix}
\usepackage[super,sort&compress]{natbib}
\usepackage{algorithm2e}
\makeatletter 
\renewcommand\@biblabel[1]{#1} 
\makeatother

\DeclareMathOperator{\TOPF}{TOP_F}
\DeclareMathOperator{\BOTF}{BOT_F}
\DeclareMathOperator{\TOPW}{TOP_W}
\DeclareMathOperator{\BOTW}{BOT_W}
\DeclareMathOperator{\Feed}{F}
\DeclareMathOperator{\Sidedraw}{W}
\DeclareMathOperator{\SEC}{SEC}

\begin{document}
\title{Minimum reflux calculation for multicomponent distillation in multi-feed, multi-product columns: Algorithms and examples}
\author[1,3]{Zheyu Jiang$^\star$}
\author[2]{Mohit Tawarmalani$^*$}
\author[1]{Rakesh Agrawal$^\dagger$}
\affil[1]{\normalsize Davidson School of Chemical Engineering, Purdue University, West Lafayette, IN 47907}
\affil[2]{Mitch Daniels School of Business, Purdue University, West Lafayette, IN 47907}
\affil[3]{School of Chemical Engineering, Oklahoma State University, Stillwater, OK 74078 \vspace{-6ex}}

\date{}

\maketitle
\noindent Corresponding authors: 
\texttt{zjiang@okstate.edu}$^\star$,
\texttt{mtawarma@purdue.edu}$^*$,
\texttt{agrawalr@purdue.edu}$^\dagger$

\begin{abstract}
\vspace{-1em}
In this work, we present the first algorithm for identifying the minimum reboiler vapor duty requirement for a general multi-feed, multi-product (MFMP) distillation column separating ideal multicomponent mixtures. This algorithm incorporates our latest advancement in developing the first shortcut model for MFMP columns. We demonstrate the accuracy and efficiency of this algorithm through case studies. The results obtained from these case studies also provide valuable insights on optimal design of MFMP columns. Many of these insights are against the existing design guidelines and heuristics. For example, placing a colder saturated feed stream above a hotter saturated feed stream sometimes leads to higher energy requirement. Furthermore, decomposing a general MFMP column into individual simple columns may lead to incorrect estimation of the minimum reflux ratio for the MFMP column. Thus, the algorithm presented here offers a fast, accurate, and automated approach to synthesize new, energy-efficient, and cost-effective MFMP columns.

\end{abstract}

\section{Introduction}

Distillation is a ubiquitous separation technology in the chemical process industries, consuming almost 50\% of the energy used by the chemical industries and about 40\% by the refining process \cite{distillationbook}. Assuming that 50\% of the $\text{CO}_2$ equivalent release from process heating in chemical manufacturing and 40\% in petroleum refining are attributable to distillation, distillation alone would be responsible for 95 million tons of $\text{CO}_2$ release in the U.S. each year \cite{amo}. Thus, to decarbonize the U.S. manufacturing sector, it is essential to significantly reduce the energy consumption and carbon footprint of distillation process \cite{roadmap}.

\begin{figure}[!ht]
	\centering
	\includegraphics[width=0.9\textwidth]{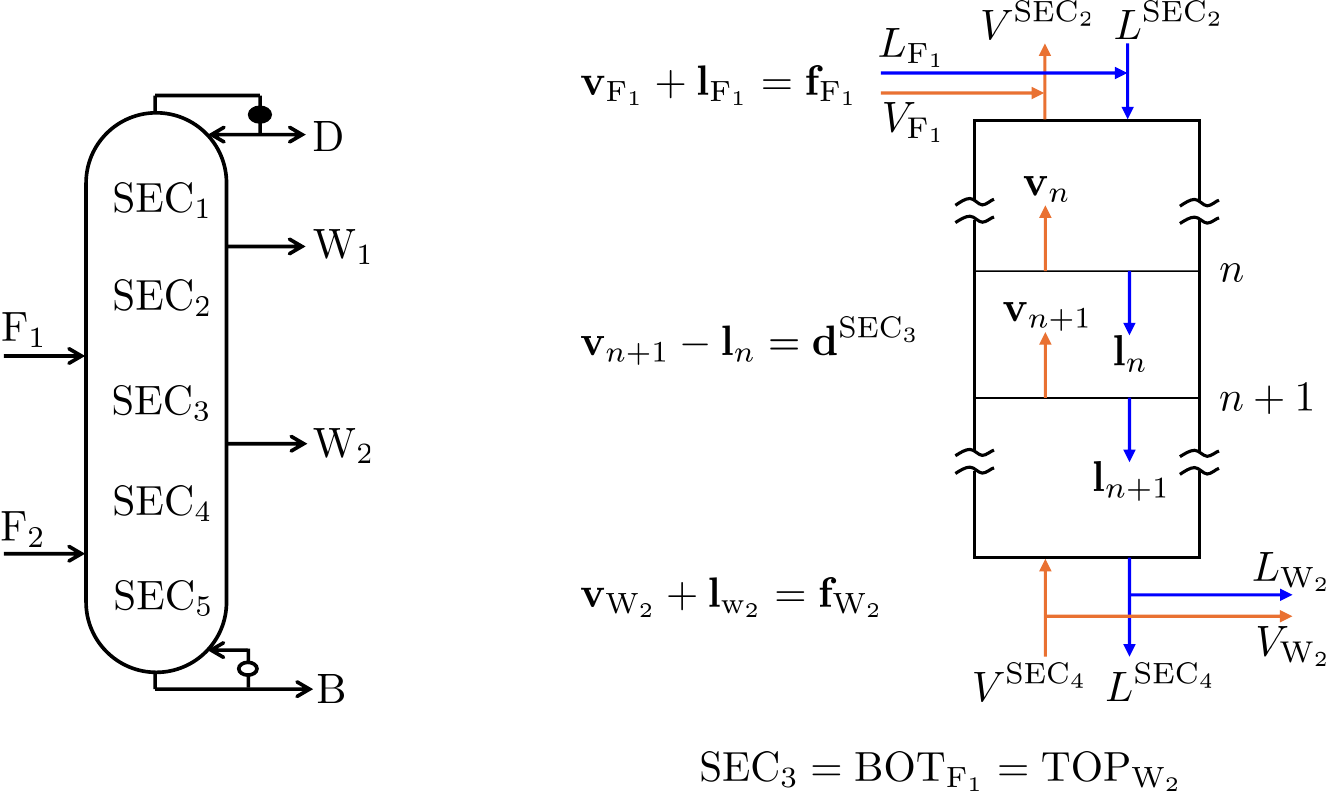}
	\caption{An example MFMP column with three feed streams and two sidedraw product streams and a detailed illustration of liquid and vapor flows within $\SEC_3$, in which the variables in bold are component flow vectors (e.g., $\mathbf{d}^{\sec_3} = (d_1^{\sec_3},  \dots, d_c^{\sec_3})$ for a $c$-component system). The column section is numbered from top (1) to bottom ($N_{\SEC}=5$). The definitions of variables and parameters used here and for the rest of this paper are summarized in Appendices \ref{appendixA} and \ref{appendixB}. We follow the convention that $\mathbf{v}_{\Feed_j}$, $\mathbf{l}_{\Feed_j}$, and $\mathbf{f}_{\Feed_j} \geq \mathbf{0}$, whereas $\mathbf{v}_{\Sidedraw_j}$, $\mathbf{l}_{\Sidedraw_j}$, and $\mathbf{f}_{\Sidedraw_j} \leq \mathbf{0}$. }
	\label{fig_mfmp}
\end{figure}

While binary mixtures can generally be separated using one distillation column, multicomponent mixtures, which are more commonly encountered in industrial separations, require a sequence of columns called a distillation configuration to achieve the desired separation. As the number of components in the feed increases, the total number of possible distillation configurations increases combinatorially \cite{shahmatrix}. Among these distillation configurations, many of them contain one or more distillation columns with multiple feed streams and/or one or more sidedraw product streams. And it is well-known that these configurations with multi-feed, multi-product (MFMP) columns (see Figure \ref{fig_mfmp} for an example) are generally more energy-efficient to operate than the so called ``sharp-split configurations'' which do not involve any MFMP column \cite{gmaneed,nallasivam}.

MFMP columns can also be derived from conventional one-feed, two-product columns in binary and multicomponent distillation by applying various process intensification techniques \cite{pireview,pireview2,misconception}, including heat pumps \cite{heatpump, akash}, double and multi-effect \cite{wankat}, intermediate reboilers and condensers \cite{intermediatereboiler}, prefractionator arragement \cite{prefractionation}, feed preconditioning \cite{feedconditioning}, heat and mass integration \cite{hmp}, and so on. Compared to the original conventional columns, these news MFMP columns not only require significantly less energy from a first-law of thermodynamics perspective, but also have much higher thermodynamic efficiency from a second-law perspective \cite{krishna}, making them more attractive than alternative technologies (e.g., membranes) for a variety of industrial separations \cite{misconception,jose}. Furthermore, when heat pumps are used in conjunction with other techniques above, the resulting MFMP columns can now be flexibly powered by alternative energy sources (e.g., renewable electricity such as solar and wind). Thus, MFMP columns are becoming increasingly important in the context of industrial decarbonization and net-zero economy, as they can revamp conventional steam-driven distillation systems whose energy primarily comes from fossil fuel combustion \cite{amo}.

The minimum reflux ratio of a distillation column is closely related to its energy consumption, capital cost, and operational limit \cite{gilliland2,doherty}, hence it is a key parameter in distillation design and operation. The naive approach of determining a column's minimum reflux ratio involves performing exhaustive sensitivity analysis using process simulators, which is a tedious task that often faces convergence issues. As a result, a fast and accurate algorithmic approach to calculate the actual minimum reflux condition of a general MFMP column is critical for designing new, energy-efficient, and cost-effective multicomponent distillation systems. Ideally, such a method should also have a simple mathematical formulation that can be easily incorporated in a (global) optimization framework for fast and accurate identification of attractive configurations from an enormous configuration search space.

Over the past decades, a number of algorithmic methods have been proposed to determine the minimum reflux ratio of a general MFMP column accurately and efficiently. A comprehensive review of these methods can be found in the first article of this series \cite{jiang}. However, these methods either rely on several simplifying assumptions, some of which turn out to be too restrictive or even incorrect as we will later demonstrate, or they require rigorous tray-by-tray calculations which are computationally expensive to perform and thus impractical to be implemented for solving complex MFMP columns. To fill the gap between existing methods and what practitioners anticipate, in our previous work \cite{jiang}, we develop the first shortcut mathematical model to analytically determine the minimum reflux ratio of any general MFMP column entirely based on the assumptions of ideal vapor-liquid equilibrium, constant relative volatility, and constant molar overflow. Our shortcut model is fully generalized as it works for any MFMP column with no particular requirement on feed and/or sidedraw arrangement or product composition specification. Also, the proposed shortcut model does not involve any tray-by-tray calculations. Furthermore, the physical and mathematical properties associated with the shortcut model are explored, from which we successfully derive the mathematical conditions for any general MFMP column operated at minimum reflux. In addition, a relaxation of the constant molar overflow assumption was proposed recently without changing the general mathematical structure of the governing equation \cite{cht}, hence extending the applicability of our shortcut model to real multicomponent systems even further while preserving the mathematical properties and minimum reflux conditions of our shortcut model.

Continuing our previous work \cite{jiang}, in this article, we introduce an algorithmic method that incorporates the shortcut model developed earlier to efficiently and accurately determine the minimum reboiler vapor duty requirement for a general MFMP column separating a multicomponent mixture. This algorithm can either be used by itself to find the minimum reflux condition for a standalone MFMP column, or can be embedded into a global optimization framework \cite{nallasivam, nlpcost, nlpexergy, krishna} to simultaneously optimize an entire configuration consisting of one or more MFMP columns. Later, we present three case studies in comparison with rigorous Aspen Plus simulations to illustrate the accuracy and usefulness of our algorithm. Also, we show that results from these case studies could challenge some of the widely-used design heuristics and rules-of-thumb researchers and practitioners have been relying upon. Thus, our shortcut method and the minimum reflux calculation algorithm provide new perspectives on how to accurately model, design, and operate MFMP columns. 

\section{A Brief Summary of Shortcut Model for MFMP Columns}\label{sec:model}

Before we introduce the minimum reflux calculation algorithm for MFMP columns, we first present a high-level review of the shortcut model we developed earlier \cite{jiang} and some of the key results resulting from its derivation, including the mathematical conditions that dictate whether the target separation task can be achieved (with finite or inifite number of stages) in the MFMP column. We encourage readers to refer back to our previous work \cite{jiang} for detailed derivations and explanation of these results. We consider a column section, which is separated by either a feed or a product stream, as the smallest module of a MFMP column. The idea is that, by constructing a shortcut model of a column section and exploring its mathematical and physical properties, we can derive a set of algebraic constraints that must be satisfied for each and every pair of adjacent column sections to maintain connectivity of the (liquid) composition profile between any adjacent sections, hencing enforcing the feasibility of separation of the entire MFMP column. In particular, when the target separation can only be achieved by requiring an infinite number of stages (i.e., some column sections have to be pinched), then the corresponding reflux ratio is the minimum reflux ratio of the MFMP column with respect to the target separation goal.

Consider a MFMP column with $N_{\SEC}$ number of column sections separated by $N_{\Feed}$ number of feed and $N_{\Sidedraw}$ number of sidedraw streams (note that $N_{\SEC} = N_{\Feed} + N_{\Sidedraw} + 1$). Following the nomencalture used in our previous work \cite{jiang} and herein smmarized in Appendices \ref{appendixA} and \ref{appendixB}, for a $c$-component system, let $\mathcal{C} = \{1,\dots,c\}$ and $\alpha_c > \alpha_{c-1} > \dots > \alpha_1 =1 $ be the relative volatilities with respect to the least volatile component (component 1). Given the feed and product flow rate and composition specifications, we can determine the net material upward flow for component $i$ in column section $k$, namely $d^{{\SEC}_k}_i$ (see Figure \ref{fig_mfmp}). Then, for a specific section vapor flow $V^{\SEC_k}$, we can solve the following equation \cite{jiang} to obtain a total of $c$ roots, $\{\gamma^{\SEC_k}_i\}_{i \in \mathcal{C}}$:
\begin{equation} \label{eqn_characteristic} 
\sum^{c}_{i=1} \frac{\alpha_i d_i^{\SEC_k}}{\alpha_i - \gamma^{\SEC_k}} = V^{\SEC_k}.
\end{equation}

Suppose $d_c,\dots,d_l>0$, $d_{l-1},\dots,d_{h+1}=0$, and $d_{h},\dots,d_1<0$ for some $1 \leq h<l \leq c$ in a column section. In other words, more volatile components $c,\dots, l$ have net material upward flows, intermediately volatile components $l-1,\dots,h+1$ have net zero material flows, and less volatile components $h,\dots,1$ have net material downward flows. It can be verified that all $c$ roots lie in the following intervals:
\begin{equation} \label{eqn_roots}
	\begin{aligned}
		\gamma^{\SEC_k}_i \in (\alpha_i,\alpha_{i+1}) \qquad &\text{for }i\in \{1,\dots,h\} \\
		\gamma^{\SEC_k}_i = \alpha_i \qquad &\text{for }i \in \{h+1,\dots,l-1\}\\
		\gamma^{\SEC_k}_i \in (\alpha_{i-1},\alpha_i) \qquad &\text{for }i \in \{l,\dots,c\}.
	\end{aligned}
\end{equation}

\begin{figure}[!ht]
	\centering
	\includegraphics[scale=0.6]{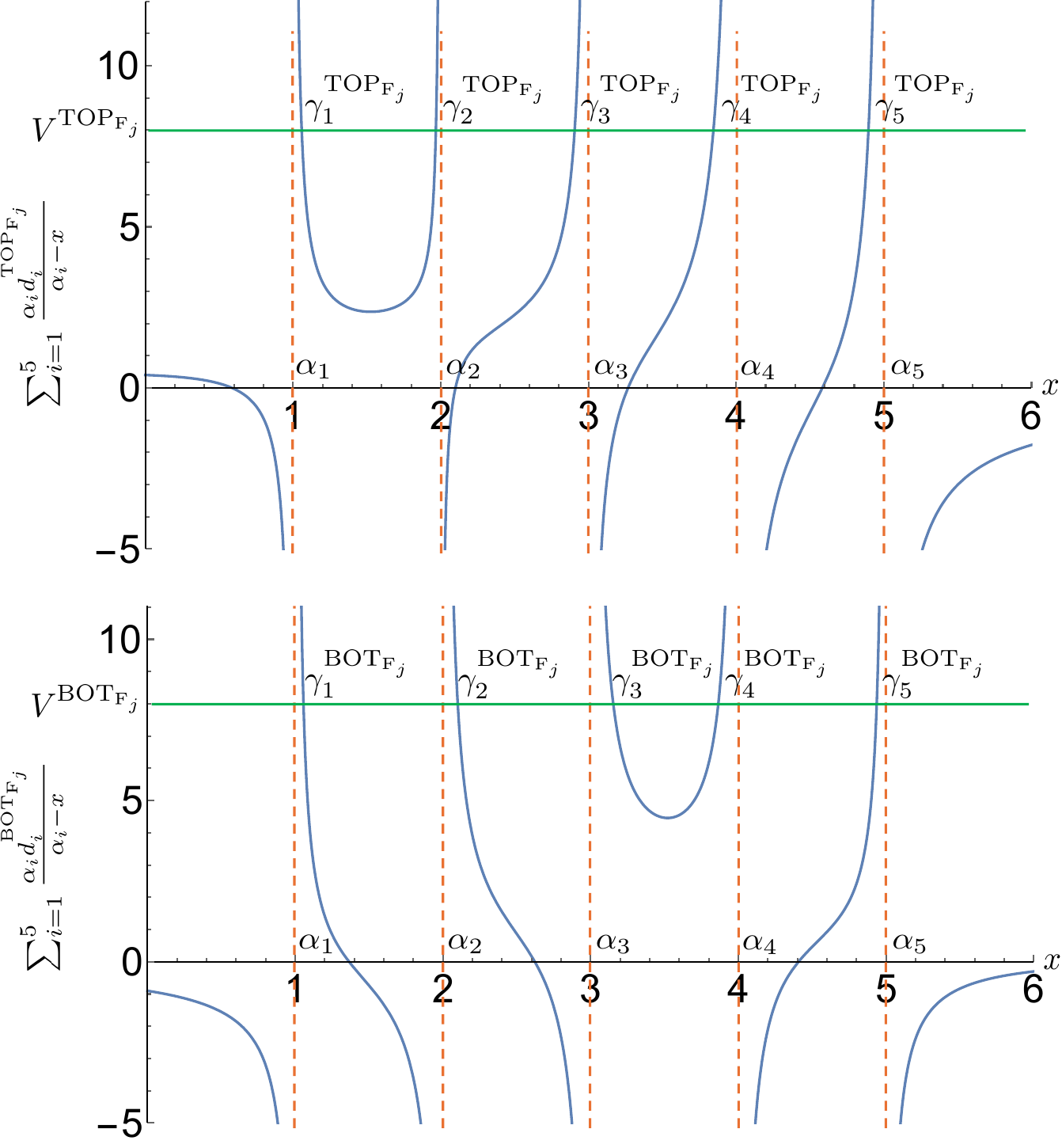}
	\caption{Roots of Equation \eqref{eqn_characteristic} for a five-component illustrative example where $(d_1, d_2, d_3, d_4, d_5) = (-0.4, 0.1, 0.2, 0.3, 0.2)$ for section $\TOPF_j$, $(d_1, d_2, d_3, d_4, d_5) = (-0.4, 0.1, 0.2, 0.3, 0.2)$ for feed $\TOPF_j$, and $(d_1, d_2, d_3, d_4, d_5) = (-0.5, -0.4, -0.3, 0.2, 0.1)$ for section $\BOTF_j$. The relative volatilities are $(\alpha_1,\alpha_2,\alpha_3,\alpha_4,\alpha_5) = (1,2,3,4,5)$. The section vapor flow $V$ is set to be 8 and $\Feed_j$ is a saturated liquid. In this case, the pinch roots $\gamma_p^{\TOPF_j} \in (\alpha_1,\alpha_2)$ and $\gamma_p^{\BOTF_j} \in (\alpha_3,\alpha_4)$.}
	\label{fig_shortcutdemo}
\end{figure}

As for the edge cases, when $l = h+1$, meaning that there are no intermediate components with net zero material flows, there are two roots in the interval $(\alpha_h,\alpha_l) = (\alpha_h,\alpha_{h+1})$ and exactly one root in each of the remaining $c-1$ relative volatility intervals (see Figure \ref{fig_shortcutdemo}). It turns out that the pinch root $\gamma^{\SEC_k}_p$ (the subscript $p$ stands for ``pinch'' and its value corresponds to the pinch index), which determines the actual pinch zone composition in $\SEC_k$, actually lies in $(\alpha_h,\alpha_l)$. That is, the pinch root lies in relative volatility interval where the sign change in $d_i$ occurs \cite{jiang}.

For the second edge case, when $h=0$, meaning that all components have non-negative net material upward flow in $\SEC_k$, we have $\gamma^{\SEC_k}_i = \alpha_i$ for $i\in \{1, \dots, l-1\}$ and $\gamma^{\SEC_k}_i \in (\alpha_{i-1}, \alpha_i)$ for $i\in \{l, \dots, c\}$. In this case, the pinch root $\gamma_p^{\SEC_k} = \gamma^{\SEC_k}_l \in (\alpha_{l-1}, \alpha_l)$. Additionally, if $l=1$, meaning that all components have net material upward flow in the column section, we have $\gamma^{\SEC_k}_i \in (\alpha_{i-1}, \alpha_i)$ for $i\in \mathcal{C}$, where $\alpha_0$ is defined as 0. And the pinch root $\gamma_p^{\SEC_k} = \gamma^{\SEC_k}_1 \in (\alpha_0, \alpha_1)$.

Lastly, when $l = c+1$, meaning that all components have non-positive net material upward flow in $\SEC_k$, we have $\gamma^{\SEC_k}_i \in (\alpha_i, \alpha_{i+1})$ for $i\in \{1, \dots, h\}$ and $\gamma^{\SEC_k}_i = \alpha_i $ for $i\in \{h+1, \dots, c\}$. In this case, the pinch root $\gamma_p^{\SEC_k} = \gamma^{\SEC_k}_h \in (\alpha_{h}, \alpha_{h+1})$. Additionally, if $h=c$, meaning that all components have net material downward flow in the column section, we have $\gamma^{\SEC_k}_i \in (\alpha_{i}, \alpha_{i+1})$ for $i\in \mathcal{C}$. Here, we denote $\alpha_{c+1} = \alpha_c + \delta $, where $\delta$ is set to be a sufficient large number. And the pinch root $\gamma_p^{\SEC_k} = \gamma^{\SEC_k}_c \in (\alpha_c, \alpha_{c+1})$.

Now that we have reviewed the key results of our shortcut model, in the next section, we will derive the algorithmic formulation to determine the minimum reflux ratio or minimum reboiler vapor duty for a general MFMP column. 

\section{Minimum Reflux Condition Formulation for MFMP Columns}\label{sec:formulation}

Recall that for $c$-component system, the domain of $\gamma_i^{\SEC_k}$ roots to Equation \eqref{eqn_characteristic} can be split into $c+1$ distinct intervals: $(0,\alpha_1)$, $(\alpha_1,\alpha_2)$, $\dots$, $(\alpha_{c-1},\alpha_c)$, and $(\alpha_c,\alpha_c+\delta)$, where $\delta$ is a sufficiently large positive number. The pinch root $\gamma_p^{\SEC_k}$ for $\SEC_k$ must lie in one of these $c+1$ intervals. Thus, we may define a set of binary variables $\{\mu_i^{\SEC_k} \in \{0,1\}\}_{i=1}^{c+1}$, where $\mu_i^{\SEC_k} = 1$ when the pinch root $\gamma_p^{\SEC_k} \in (\alpha_{i-1},\alpha_{i})$, and is 0 otherwise. Here, we denote $\alpha_0 = 0$. This way, the pinch root must satisfy the following constraints:
\begin{equation} \label{eqn_newbounds}
	\begin{aligned}
		&\sum^{c+1}_{i=1}\alpha_{i-1}\mu_i^{\SEC_k} \leq \gamma_p^{\SEC_k} \leq \sum^{c+1}_{i=1}\alpha_{i}\mu_i^{\SEC_k} \\
		&\sum^{c+1}_{i=1} \mu_i^{\SEC_k} =1
	\end{aligned}
	\qquad \forall k=1,\dots, N_{\SEC}.
\end{equation}

When two adjacent column sections are separated by a feed stream ${\Feed}_j$ ($j = 1,\dots,N_{\Feed}$), we denote the section above ${\Feed}_j$ as $\TOPF_{j}$ and the one below as $\BOTF_j$ (see Figure \ref{fig_mfmp}). Since $d_i^{\TOPF_j} - d_i^{\BOTF_j} = f_{i,\Feed_j} \geq 0$, we have $d_i^{\TOPF_j} \geq d_i^{\BOTF_j}$ for $i \in \mathcal{C}$, indicating that the pinch index for $\TOPF_j$ would be at most the same as the pinch index for $\BOTF_j$, hence satisfying $p^{\TOPF_j} \leq p^{\BOTF_j}$, or:
\begin{equation} \label{eqn_newinequality}
	\sum^{c+1}_{i=1}i\mu_{i,\BOTF_j} \geq \sum^{c+1}_{i=1}i\mu_{i,\TOPF_{j}} \quad \forall j=1,\dots,N_{\Feed}.
\end{equation}

We define an index set $\mathcal{I}_{\Feed_j}$ storing all indices of intervals ranging from $\max\{2,\sum^{c+1}_{i=1}i\mu_{i,\TOPF_{j}}\}$ to $\min\{c, \sum^{c+1}_{i=1}i\mu_{i,\BOTF_j}\}$ (i.e., considering intervals within $\alpha_1$ and $\alpha_c$ and excluding the two intervals $(\alpha_0, \alpha_1)$ and $(\alpha_c,\alpha_{c+1})$). To better characterize $\mathcal{I}_{\Feed_j}$, we define a new set of binary variables $\{K_i^{\SEC_k} \in \{0,1\} \}_{i=1}^{c+1}$ for column section $k$ where:
\begin{equation} \label{eqn_K}
	K^{\SEC_k}_{i} = \sum_{m=1}^{i} \mu_m^{\SEC_k} \qquad \forall i = 1,\dots, c+1;\; k=1,\dots,N_{\SEC}.
\end{equation}

Clearly, $K_{i}^{\SEC_k} = 0$ if and only if $\mu^{\SEC_k}_1,\dots,\mu^{\SEC_k}_{i}$ are all equal to 0. And $K_{i}^{\SEC_k}$ changes from 0 to 1 at index $i$ where $\mu_{i}^{\SEC_k} = 1$ (i.e., $\gamma_p^{\SEC_k} \in (\alpha_{i-1}, \alpha_{i})$) and then stays at 1 for indices greater than $i$. Knowing this, it can be verified that $\mathcal{I}_{\Feed_j}$ can be equivalently expressed as:
\begin{equation} \label{eqn_Ifeed}
	\mathcal{I}_{\Feed_j} = \{i \in \mathcal{C} | K_{i}^{\TOPF_j} - K_{i-1}^{\BOTF_j} = 1\} \qquad \forall j=1,\dots,N_{\Feed}.
\end{equation}

\begin{figure}[!ht]
	\centering
	\includegraphics[width=0.9\textwidth]{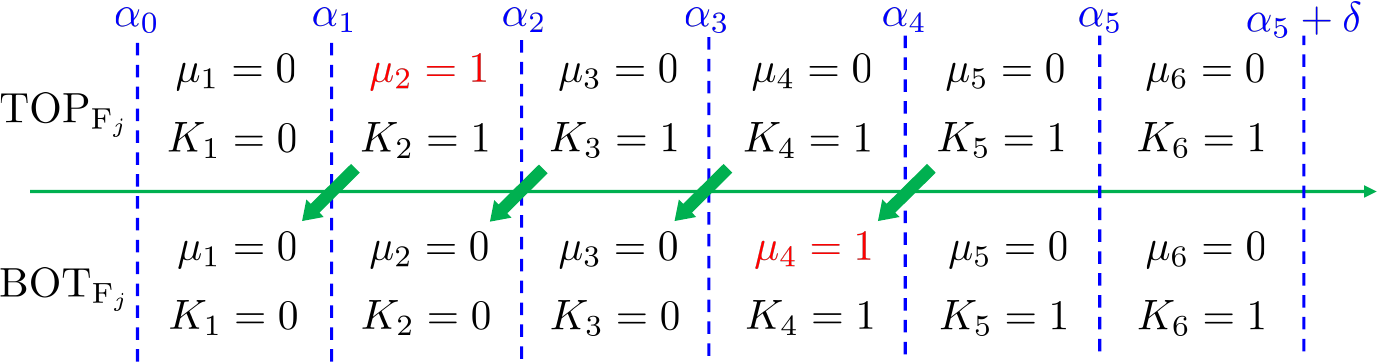}
	\caption{The relationship between $\mu_i$ and $K_i$ variables for the example illustrated in Figure \ref{fig_shortcutdemo}. The green arrows show how the (binary) coefficients in Equation \eqref{eqn_feasibility_K_feed} are constructed. For this example, $K_2^{\TOPF_j} - K_1^{\BOTF_j} = K_3^{\TOPF_j} - K_2^{\BOTF_j} = K_4^{\TOPF_j} - K_3^{\BOTF_j} = 1$. Therefore, $\mathcal{I}_{\Feed_j} = \{2,3,4\}$ according to Equation \eqref{eqn_Ifeed}.} \label{fig_K}
\end{figure}

For example, consider the same five-component system whose root profiles for $\TOPF_j$ and $\BOTF_j$ are shown in Figure \ref{fig_shortcutdemo}. Correspondingly, the relationship between $mu_i$ and $K_i$ variables for this illustrative example is shown in Figure \ref{fig_K}. Following Equation \eqref{eqn_Ifeed}, we have $\mathcal{I}_{\Feed_j} = \{2,3,4\}$. 

With this, one of the key results obtained in our previous work \cite{jiang} is that, the feasibility of the target separation requires the following constraint to the satisfies in sections $\TOPF_j$ and $\BOTF_j$ for every $i \in \mathcal{I}_{\Feed_j}$:
\begin{equation} \label{eqn_feasibility}
	\gamma_{i}^{\TOPF_j} \geq \rho_{i-1,\Feed_j} \geq \gamma_{i-1}^{\BOTF_j} \qquad \forall i \in \mathcal{I}_{\Feed_j}; \; j=1,\dots,N_{\Feed}
\end{equation}
where $\{\rho_{i-1,\Feed_j}\}_{i \in \mathcal{I}_{\Feed_j}}$ satisfy:
\begin{equation} \label{eqn_feed}
    \sum_{m=1}^c \frac{\alpha_m l_{m,\Feed_j}}{\alpha_m-\rho_{i-1,\Feed_j}}=0 \qquad \text{or} \qquad  \sum_{m=1}^c \frac{\alpha_m f_{m,\Feed_j}}{\alpha_m-\rho_{i-1,\Feed_j}} = V_{\Feed_j} \qquad j=1,\dots,N_{\Feed},
\end{equation}
where $\rho_{i-1,\Feed_j} \in (\alpha_{i-1},\alpha_{i})$. Here, $l_{m,\Feed_j} \geq 0$, $f_{m,\Feed_j} \geq 0$, and $V_{\Feed_j} \geq 0$ correspond to the flow rate of component $m$ in the liquid portion of $\Feed_j$, the feed flow rate of component $m$, and the total vapor flow rate of the feed, respectively. When $\Feed_j$ is in saturated vapor state, then $l_{m,\Feed_j}$ represents the hypothetical liquid feed flow that is in thermodynamic equilibrium with the vapor feed \cite{jiang}: 
\begin{equation*} 
	l_{m,\Feed_j} = \frac{\frac{v_{m,\Feed_j}}{\alpha_m}}{\sum_{k=1}^c \frac{v_{k,\Feed_j}}{\alpha_k}} \qquad \forall m\in \mathcal{C}.
\end{equation*}

To implement Equation \eqref{eqn_feasibility} algorithmically, we leverage the fact that $K_{i}^{\TOPF_j} - K_{i-1}^{\BOTF_j}$ is itself a binary variable indicating whether index $i\in \mathcal{I}_{\Feed_j}$ or not. Thus, Equation \eqref{eqn_feasibility} can be rewritten as:
\begin{equation} \label{eqn_feasibility_K_feed}
	\begin{aligned}
		& (K_{i}^{\TOPF_j} - K_{i-1}^{\BOTF_j})(\gamma_i^{\TOPF_j} - \rho_{i-1,\Feed_j}) \geq 0\\
		& (K_{i}^{\TOPF_j} - K_{i-1}^{\BOTF_j})(\rho_{i-1,\Feed_j} - \gamma_{i-1}^{\BOTF_j}) \geq 0
	\end{aligned}
	\qquad \forall i \in \mathcal{C}\backslash\{1\};\; j=1,\dots, N_{\Feed}.
\end{equation}

When two adjacent column sections are connected by a sidedraw stream $\Sidedraw_j$ ($j=1,\dots,N_{\Sidedraw}$), we denote the section above $\Sidedraw_j$ as $\TOPW_j$ and the one below as $\BOTW_j$. Since $d_i^{\TOPW_j} - d_i^{\BOTW_j} = f_{i,\Sidedraw_j} \leq 0$, we have $d_i^{\TOPW_j} \leq d_i^{\BOTW_j}$ for $i \in \mathcal{C}$, indicating that the pinch index for section $\TOPW_j$ would be at least the same as the pinch index for section $\BOTW_j$, hence satisfying $p^{\TOPW_j}\geq p^{\BOTW_j}$. Thus, we have:
\begin{equation} \label{eqn_newinequality2}
	\sum^{c+1}_{i=1}i\mu_{i,\BOTW_j} \leq \sum^{c+1}_{i=1}i\mu_{i,\TOPW_j} \quad \forall j=1,\dots,N_{\Sidedraw}.
\end{equation}

Similar to $\mathcal{I}_{\Feed_j}$ for feed stream $\Feed_j$, we define an index set $\mathcal{I}_{\Sidedraw_j}$ storing all indices of intervals ranging from $\max\{2,\sum^{c+1}_{i=1}i\mu_{i,\BOTW_{j}}\}$ to $\min\{c, \sum^{c+1}_{i=1}i\mu_{i,\TOPW_j}\}$. Furthermore, $\mathcal{I}_{\Sidedraw_j}$ can also be redefined as:
\begin{equation} \label{eqn_ISidedraw}
	\mathcal{I}_{\Sidedraw_j} = \{i \in \mathcal{C} | K_{i}^{\BOTW_j} - K_{i-1}^{\TOPW_j} = 1\} \qquad \forall j=1,\dots,N_{\Sidedraw}.
\end{equation}

The feasibility of separation requires the following constraint to hold for $\TOPW_j$ and $\BOTW_j$ for every $i \in \mathcal{I}_{\Sidedraw_j}$ \cite{jiang}:
\begin{equation} \label{eqn_feasibility2}
	\gamma_{i-1}^{\TOPW_j} \leq \rho_{i-1,\Sidedraw_j} \leq \gamma_{i}^{\BOTW_j} \qquad \forall i \in \mathcal{I}_{\Sidedraw_j};\; j=1,\dots,N_{\Sidedraw},
\end{equation}
where $\{\rho_{i-1,\Sidedraw_j}\}_{i \in \mathcal{I}_{\Sidedraw_j}}$ satisfy Equation \eqref{eqn_sidedraw} below:
\begin{equation} \label{eqn_sidedraw}
	\sum_{m=1}^c \frac{\alpha_m l_{m,\Sidedraw_j}}{\alpha_m-\rho_{i-1,\Sidedraw_j}} = 0, \qquad \text{or} \qquad \sum_{m=1}^c \frac{\alpha_m f_{m,\Sidedraw_j}}{\alpha_m-\rho_{i-1,\Sidedraw_j}} = V_{\Sidedraw_j} \qquad j=1,\dots,N_{\Sidedraw},
\end{equation}
where $\rho_{i-1,\Sidedraw_j} \in (\alpha_{i-1},\alpha_{i})$. Here, $l_{m,\Sidedraw_j} \leq 0$, $f_{m,\Sidedraw_j} \leq 0$, and $V_{\Sidedraw_j} \leq 0$ correspond to the flow rate of component $m$ in the liquid portion of the sidedraw stream, the sidedraw flow rate of component $m$, and the total vapor flow rate of the sidedraw, respectively. When $\Sidedraw_j$ is in saturated vapor state, then $l_{m,\Sidedraw_j}$ represents the hypothetical liquid sidedraw flow that is in thermodynamic equilibrium with the vapor sidedraw \cite{jiang}: 
\begin{equation*} 
	l_{m,\Sidedraw_j} = \frac{\frac{v_{m,\Sidedraw_j}}{\alpha_m}}{\sum_{k=1}^c \frac{v_{k,\Sidedraw_j}}{\alpha_k}} \qquad \forall m\in \mathcal{C}.
\end{equation*}

Similar to how we reformulate Equation \eqref{eqn_feasibility} for $\Feed_j$, we can rewrite Equation \eqref{eqn_feasibility2} as:
\begin{equation}\label{eqn_feasibility_K_sidedraw}
	\begin{aligned}
		& (K_{i}^{\BOTW_j} - K_{i-1}^{\TOPW_j})(\gamma_i^{\BOTW_j} - \rho_{i-1,\Sidedraw_j}) \geq 0\\
		& (K_{i}^{\BOTW_j} - K_{i-1}^{\TOPW_j})(\rho_{i-1,\Sidedraw_j} - \gamma_{i-1}^{\TOPW_j}) \geq 0
	\end{aligned}
	\qquad \forall i \in \mathcal{C}\backslash\{1\};\; j=1,\dots, N_{\Sidedraw}.
\end{equation}

Furthermore, as illustrated in Figure \ref{fig_mfmp}, a unique feature about a sidedraw stream is that the sidedraw's liquid (resp. vapor) composition must lie on the liquid (resp. vapor) composition profile, whereas a feed stream's liquid (resp. vapor) composition may or may not lie on the liquid (resp. vapor) composition profile. The condition that sidedraw composition must belong to the composition profile yields to the following set of constraints:
\begin{equation} \label{eqn_feasibility4}
	\begin{split}
		& K_i^{\TOPW_j}(\gamma_{i}^{\TOPW_j} - \rho_{i-1,\Sidedraw_j}) \geq 0  \qquad \forall i\in \mathcal{C}\backslash \{1\};\\
		& K_i^{\BOTW_j}(\gamma_i^{\BOTW_j} - \rho_{i-1,\Sidedraw_j}) \geq 0  \qquad \forall i\in \mathcal{C}\backslash \{1\};\\
		& (1-K_i^{\TOPW_j})(\gamma_{i}^{\TOPW_j} - \rho_{i,\Sidedraw_j}) \leq 0  \qquad \forall i \in \mathcal{C};\\
		& (1-K_i^{\BOTW_j})(\gamma_{i}^{\BOTW_j} - \rho_{i,\Sidedraw_j}) \leq 0  \qquad \forall i \in \mathcal{C};\\
	\end{split}
	\qquad \forall j=1,\dots,N_{\Sidedraw}.
\end{equation}

\section{Implementation of Minimum Reflux Calculation Algorithm}\label{sec:implementation}

When implementing the algorithm developed in Section \ref{sec:formulation}, there are two approaches to consider. The first approach is to implement Equations \eqref{eqn_newbounds}, \eqref{eqn_newinequality}, \eqref{eqn_K}, \eqref{eqn_feasibility_K_feed}, \eqref{eqn_feasibility_K_sidedraw} and \eqref{eqn_feasibility4} in an optimization framework as constraints, along with the mathematical formulations of the shortcut model developed in our earlier work \cite{jiang}, which includes Equations \eqref{eqn_characteristic}, \eqref{eqn_feed}, \eqref{eqn_sidedraw}, and more. The resulting formulation is a mixed-integer nonlinear program (MINLP), which can be solved to global optimality using global solvers such as BARON \cite{baron}. We pursue this approach when only the purity or recovery of the key components in product streams are specified. In this case, the MINLP will determine the optimal distribution of other components in the product streams such that the reflux ratio or reboiler vapor duty requirement is minimized. To illustrate how this approach works, in Section \ref{sec:4compcase}, we present this formulation for an quaternary separation example in a two-feed, one-sidedraw column.

\begin{algorithm}
	\SetAlgoLined
	\caption{\texttt{Vrebmin}: Algorithm for determining the minimum reboiler vapor duty requirement of a MFMP column knowing the flow rates and compositions of feed and product streams. The minimum reflux ratio $R_{\min}$ can be readily calculated from vapor balances once $V_{\text{reb},\min}$ is obtained, since all feed and product flow rates are specified.} \label{algo1}
	\SetKwData{Left}{left}\SetKwData{This}{this}\SetKwData{Up}{up}
	\SetKwFunction{Union}{Union}\SetKwFunction{FindCompress}{FindCompress}
	\SetKwInOut{Input}{input}\SetKwInOut{Output}{output}
	\Input{$\mathcal{C},\, N_{\Feed}, \, N_{\Sidedraw},\, N_{\SEC}, \, \{f_{i,\Feed_j}\}_{j=1}^{N_{\Feed}},\, \{l_{i,\Feed_j}\}_{j=1}^{N_{\Feed}},\, \{f_{i,\Sidedraw_j}\}_{j=1}^{N_{\Sidedraw}}, \, \{l_{i,\Sidedraw_j}\}_{j=1}^{N_{\Sidedraw}}, \,d_i^{\SEC_1}$ for every component $i \in \mathcal{C}$}
	\Output{minimum reboiler vapor duty $V_{\text{reb},\min}$}
	\textbf{initialize:} An empty list $\{V_{\text{reb}}\}$ storing candidate minimum reboiler vapor duty values
	\BlankLine
	\Begin{
	Calculate $\{d_i^{\SEC_k}\}_{i\in\mathcal{C};\; k=2,\dots,N_{\SEC}}$ from inter-column section material balances (Equation (29) of Jiang et al. \cite{jiang})\;
	Determine pinch root $\{\gamma_p^{\SEC_k}\}_{k=1,\dots,N_{\SEC}}$ locations from Equation \eqref{eqn_roots} and other edge cases, then obtain $\{\mu_i^{\SEC_k}\}_{i\in \mathcal{C},k=1,\dots,N_{\SEC}}$ and $\{K_i^{\SEC_k}\}_{i\in \mathcal{C},k=1,\dots,N_{\SEC}}$\;
	Determine index sets $\{\mathcal{I}_{\Feed_j}\}_{j=1}^{N_{\Feed}}$ and $\{\mathcal{I}_{\Sidedraw_j}\}_{j=1}^{N_{\Sidedraw}}$ from Equations \eqref{eqn_Ifeed} and \eqref{eqn_ISidedraw}\;
	Solve Equations \eqref{eqn_feed} and \eqref{eqn_sidedraw} to obtain $\{\rho_{i-1,\Feed_j}\}_{i \in \mathcal{C}\backslash\{1\}; \; j=1,\dots,N_{\Feed}}$ and $\{\rho_{i-1,\Sidedraw_j}\}_{i \in \mathcal{C}\backslash\{1\}; \; j=1,\dots,N_{\Sidedraw}}$ roots, respectively\;
	\For{$j \leftarrow 1$ \KwTo $N_{\Sidedraw}$}{
		Add $V_{\text{reb},\Sidedraw_j} = \texttt{sidedrawFeasible}(j,\{K_i^{\TOPW_j}\}_i)$ into the list $\{V_{\text{reb}}\}$ and continue\;
		\leIf{$\mathcal{I}_{\Sidedraw_j} = \emptyset$}{Skip and go to the next $j$}{Continue}
		\For{$i\in\mathcal{I}_{\Sidedraw_j}$}{
			Substitute $\gamma_{i-1}^{\TOPW_j} \leftarrow \rho_{i-1,\Sidedraw_j}$ into Equation \eqref{eqn_characteristic} to obtain $V^{\TOPW_j}$\;
			Add $V_{\text{reb},\Sidedraw_j} = \texttt{getVreb}(\TOPW_j, V^{\TOPW_j})$ into the list $\{V_{\text{reb}}\}$\;
		}
	}
	\For{$j \leftarrow 1$ \KwTo $N_{\Feed}$}{
		\leIf{$\mathcal{I}_{\Feed_j} = \emptyset$}{Skip and go to the next $j$}{Continue}
		\For{$i\in\mathcal{I}_{\Feed_j}$}{
			Substitute $\gamma_{i}^{\TOPF_j} \leftarrow \rho_{i-1,\Feed_j}$ into Equation \eqref{eqn_characteristic} to obtain $V^{\TOPF_j}$\;
			Add $V_{\text{reb},\Feed_j} = \texttt{getVreb}(\TOPF_j, V^{\TOPF_j})$ into the list $\{V_{\text{reb}}\}$\;
		}
	}
	$V_{\text{reb},\min} = \min\{V_{\text{reb}}\}$
	}
\end{algorithm}

\begin{algorithm}
	\SetAlgoLined
	\caption{\texttt{getVreb}: Algorithm for checking feasibility of separation and returning the candidate reboiler vapor duty value.} \label{algo2}
	\SetKwData{Left}{left}\SetKwData{This}{this}\SetKwData{Up}{up}
	\SetKwFunction{Union}{Union}\SetKwFunction{FindCompress}{FindCompress}
	\SetKwInOut{Input}{input}\SetKwInOut{Output}{output}
	\Input{column section $s$ and its section vapor flow $V^{\SEC_s}$}
	\Output{candidate reboiler vapor duty value}
	\BlankLine
	\Begin{
		\For{$k \leftarrow s - 1$ \KwTo $1$} {
			Calculate $V^{\SEC_k}$ from $V^{\SEC_{k+1}}$ via vapor balances\;
			Determine $\{\gamma_r^{\SEC_k}\}_{r \in \mathcal{C}}$ from Equation \eqref{eqn_characteristic} using $V^{\SEC_k}$\;
			\leIf{the feasibility constraints (Equations \eqref{eqn_feasibility}, \eqref{eqn_feasibility2}, \eqref{eqn_feasibility4}) associated with $\SEC_k$ are satisfied}{Continue}{\Return{\texttt{null}}}
			}
			\For{$k \leftarrow s + 1$ \KwTo $N_{\SEC}$} {
			Calculate $V^{\SEC_k}$ from $V^{\SEC_{k-1}}$ via vapor balances\;
			Determine $\{\gamma_r^{\SEC_k}\}_{r \in \mathcal{C}}$ from Equation \eqref{eqn_characteristic} using $V^{\SEC_k}$\;
			\leIf{the feasibility constraints (Equations \eqref{eqn_feasibility}, \eqref{eqn_feasibility2}, \eqref{eqn_feasibility4}) associated with $\SEC_k$ are satisfied}{Continue}{\Return{\texttt{null}}}
			}
		\Return{$V^{\SEC_{N_{\SEC}}}$}.
	}
\end{algorithm}

\begin{algorithm}
	\SetAlgoLined
	\caption{\texttt{sidedrawFeasible}: Algorithm for returning candidate reboiler vapor duty value assuming that the minimum reflux is ``controlled'' by a sidedraw.} \label{algo3}
	\SetKwData{Left}{left}\SetKwData{This}{this}\SetKwData{Up}{up}
	\SetKwFunction{Union}{Union}\SetKwFunction{FindCompress}{FindCompress}
	\SetKwInOut{Input}{input}\SetKwInOut{Output}{output}
	\Input{sidedraw stream index $j$ and index set $\{K_i^{\TOPW_j}\}_i$}
	\Output{candidate reboiler vapor duty value}
	\textbf{initialize:} An empty list $\{V_{\text{reb},\Sidedraw_j}\}$ storing candidate minimum reboiler vapor duty values
	\BlankLine
	\Begin{
		\For{$m \leftarrow 1$ \KwTo $c$}{
			\uIf{$K_m^{\TOPW_j} = 0$}{Substitute $\gamma_m^{\TOPW_j} \leftarrow \rho_{m,\Sidedraw_j}$ into Equation \eqref{eqn_characteristic} to obtain $V^{\TOPW_j}$; Add $V_{\text{reb},m} = \texttt{getVreb}(\TOPW_j, V^{\TOPW_j})$ into the list $\{V_{\text{reb},\Sidedraw_j}\}$\;}
			\uElse{Substitute $\gamma_m^{\TOPW_j} \leftarrow \rho_{m-1,\Sidedraw_j}$ into Equation \eqref{eqn_characteristic} to obtain $V^{\TOPW_j}$; Add $V_{\text{reb},m} = \texttt{getVreb}(\TOPW_j, V^{\TOPW_j})$ into the list $\{V_{\text{reb},\Sidedraw_j}\}$\;}
		}
		\Return{$V^{\SEC_{N_{\SEC}}} \leftarrow \min \{V_{\text{reb},\Sidedraw_j}\}$}.
	}
\end{algorithm}

The second approach deals with many practical applications in which the product distributions of the MFMP column have already been adequately specified. In this case, the search for minimum reflux ratio of the MFMP column becomes a fully algorithmic procedure that does not require solving an optimization problem. This is because the net material upward flows $\{d_i^{\SEC_k}\}_{i\in \mathcal{C},k=1,\dots,N_{\SEC}}$ can be readily obtained from mass balances, making the determination of pinch root location (hence $\mathcal{I}_{\Feed}$ and $\mathcal{I}_{\Sidedraw}$ sets) completely deterministic for every feed and sidedraw stream. Therefore, we can run a simple algorithmic procedure, as shown in Algorithms \ref{algo1} through \ref{algo3}, to identify the true minimum reboiler vapor duty requirement (or equivalently the minimum reflux ratio since all product flows are fixed). Specifically, as discussed in detail in our previous work \cite{jiang}, at minimum reflux condition, one of the feed or sidedraw streams essentially ``controls'' the separation. Accordingly, while the feasibility constraints (Equation \eqref{eqn_feasibility} or \eqref{eqn_feasibility2}) associated with feed and/or sidedraw streams continue to hold, the feasibility constraints associated with the controlling feed or sidedraw stream will become binding (i.e., satisfied as equalities). Thus, the idea behind Algorithms \ref{algo1} through \ref{algo3} is to scrutinize all feed and sidedraw streams, assuming that each of them might be ``controlling'' the separation at minimum reflux, and determine whether feasibility constraints are met for the rest of feed and sidedraw streams. Overall, the true reboiler vapor duty (resp. minimum minimum reflux ratio) corresponds to the lowest reboiler vapor duty (resp. lowest reflux ratio) at which all feasibility constraints are satisfied.

\section{Case Studies}\label{sec:examples}

In this section, we examine a few ternary and quanterary separation examples that will illustrate the accuracy and effectiveness of our minimum reflux calculation methods while providing valuable insights of the minimum reflux behavior of a MFMP column for the first time.

\subsection{Example 1: Two-Feed Distillation Column}

\begin{figure}[ht]
\centering
\includegraphics[scale=1.0]{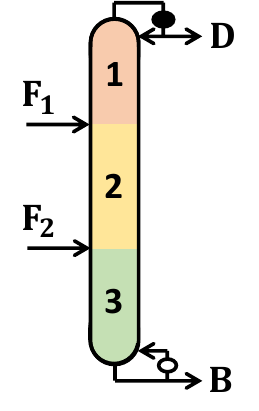}
\vspace{-1em}
\caption{A two-feed column with no sidedraw product stream.}
\label{fig_case1}
\end{figure}

In the first example, we examine a two-feed distillation column shown in Figure \ref{fig_case1} separating a ternary mixture of n-hexane (Component 3), n-heptane (Component 2), and n-octane (Component 1). Two-feed columns are common in extractive distillation applications. Furthermore, as recently discovered by Madenoor Ramapriya et al. \cite{gauthamtwofeed}, a large energy saving can potentially be realized when two feed streams are introduced at two different locations of the column compared to pre-mixing them to form a single feed stream.

The relative volatility of each component with respect to n-octane at atmospheric pressure are estimated from Aspen Plus to be $(\alpha_3,\alpha_2,\alpha_1)=(5.1168, 2.25, 1)$. To establish a common basis for comparison, we ensure constant relative volatility and constant molar overflow assumptions by appropriately modifying the property parameters in Aspen Plus listed under PLXANT and DHVLDP \cite{shah}. The IDEAL thermodynamic package is used. This column produces a distillate product with a total flow rate of 52.476 mol/s containing 95 mol\% of n-hexane, 5 mol\% of n-heptane, and negligible amount of n-octane. Thus, bottoms product has a flow rate of 147.524 mol/s containing 0.1 mol\% of n-hexane, 45.671 mol\% of n-heptane, and 54.229 mol\% of n-octane. 

We consider two scenarios in Example 1. In the first scenario, the upper feed $\Feed_1$ in the MFMP column is a saturated liquid stream containing 30 mol/s of n-hexane, 60 mol/s of n-heptane, and 10 mol/s of n-octane. The lower feed $\Feed_2$ is also a saturated liquid stream but with 20 mol/s of n-hexane, 10 mol/s of n-heptane, and 70 mol/s of n-octane. Clearly, $\Feed_2$ is less volatile (i.e., ``heavier'') than $\Feed_1$ and thus has a higher bubble point temperature. Since the feed and product specifications are given, we determine that $\mathcal{I}_{\Feed_1} = \{2,3\}$ and $\mathcal{I}_{\Feed_2} = \{3\}$ based on our earlier discussion in Section \ref{sec:model}. Substituting these index sets to Algorithms \ref{algo1} through \ref{algo3}, we obtain that the minimum reflux ratio $R_{\min}=2.162$ and the corresponding minimum reboiler vapor duty is $V_{\text{reb},\min} = 165.95$ mol/s. The minimum reflux condition occurs when the upper feed $\Feed_1$ ``controls'' the separation, in which Equation \eqref{eqn_feasibility} associated with $i = 3$ in $\mathcal{I}_{\Feed_1}$ becomes the binding constraint ($\gamma_3^{\TOPF_1} = \gamma_2^{\BOTF_1} = \rho_{2,\Feed_2}$).

\begin{figure}[!ht]
	\centering
	\includegraphics[scale=0.7]{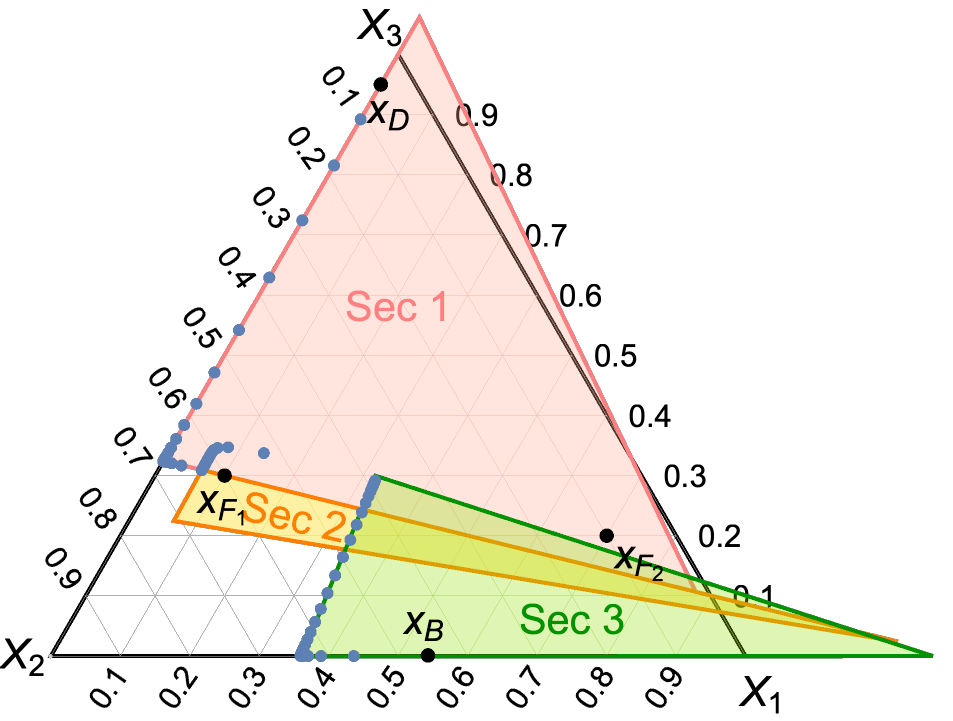}
	\caption{The pinch simplices at the minimum reflux condition obtained using Algorithms \ref{algo1} through \ref{algo3}. Hereafter, $X_1,\, X_2,\, X_3$ represent pure n-octane, n-heptane, and n-hexane, respectively. The colors of the pinch simplices match with those in Figure \ref{fig_case1}. The blue dots are the actual liquid composition profile of this two-feed column simulated in Aspen Plus as a RadFrac column. By setting up appropriate \texttt{Design Specs} in Aspen Plus to simulate the MFMP containing 150 equilibrium stages, we obtain a minimum reflux ratio of $R_{\min}=2.145$ from Aspen Plus. The exact pinches compositions in $\SEC_1$ through $\SEC_3$ are $Z_2$ (associated with pinch root $\gamma_p^{\SEC_1} = \gamma_2^{\SEC_1} \in (\alpha_1, \alpha_2)$), $Z_3$ (associated with pinch root $\gamma_p^{\SEC_2} = \gamma_3^{\SEC_2} \in (\alpha_2, \alpha_3)$), and $Z_3$ ($\gamma_p^{\SEC_3} = \gamma_3^{\SEC_3} \in (\alpha_3, \alpha_3+\delta)$), respectively. Therefore, $\mu_2^{\SEC_1} = \mu_3^{\SEC_2} = \mu_4^{\SEC_3} = 1$.}
	\label{fig_case1scenario1}
\end{figure}

\begin{figure}[!ht]
	\centering
	\includegraphics[width=\textwidth]{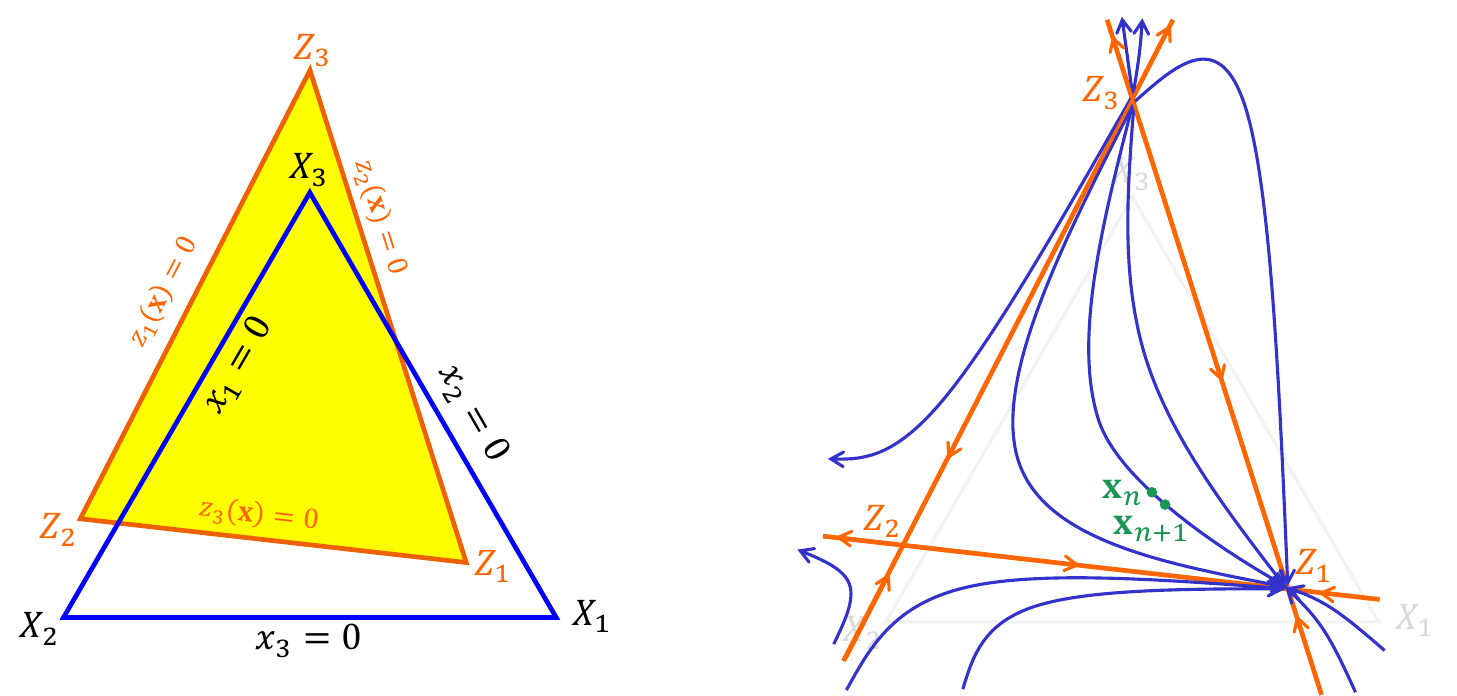}
	\caption{An illustration of pinch simplex constructed for a column section and liquid composition trajectory bundle. The pinch simplex boundary $z_i(\mathbf{x}) = 0$ is associated with the root $\gamma_i$ (see Table 1 of our previous work \cite{jiang} for explanation). And possible pinch compositions are given by vertices of the pinch simplex, $Z_i$. The arrows indicate the direction of liquid composition evolution as we move downward from the top of the column section.}
	\label{fig_simplex}
\end{figure}

For this ternary separation, we can visualize its minimum reflux condition by constructing the pinch simplicies in Figure \ref{fig_case1scenario1} based on our previous work\cite{jiang}. We observe that the pinch simplices associated with $\SEC_1$ and $\SEC_2$ share a common boundary, where $\Feed_1$ stream composition $\mathbf{x}_{\Feed_1}$ also lies. This means the two boundaries of the pinch simplices satisfy $z_3^{\TOPF_1}(\mathbf{x}_{\Feed_1}) = z_2^{\BOTF_1}(\mathbf{x}_{\Feed_1}) = 0$, which implies $\gamma_3^{\TOPF_1} = \gamma_2^{\BOTF_1} = \rho_{2,\Feed_2}$ (see Figure \ref{fig_simplex} for illustration of pinch simplex; readers may refer to our previous work\cite{jiang} for detailed explanation). In other words, the geometric interpretation of feasible separation is that the pinch simplices of any two adjacent column sections must be connected, and minimum reflux condition occurs when the pinch simplices sandwiching the controlling feed or sidedraw stream share a common face. If the reflux ratio is further reduced, these two simplices will no longer be connected, or $z_3^{\TOPF_1}(\mathbf{x}_{\Feed_1}) < 0$, and $z_2^{\BOTF_1}(\mathbf{x}_{\Feed_1}) < 0$. This indicates that $\gamma_3^{\TOPF_1} < \rho_{2,\Feed_1}$ and $\gamma_2^{\BOTF_1} > \rho_{2,\Feed_1}$, hence violating the feasibility constraint of Equation \eqref{eqn_feasibility} for feed $\Feed_1$. Therefore, $R_{\min}=2.162$ is indeed the minimum reflux ratio.

\begin{figure}[!ht]
	\centering
	\includegraphics[width=\textwidth]{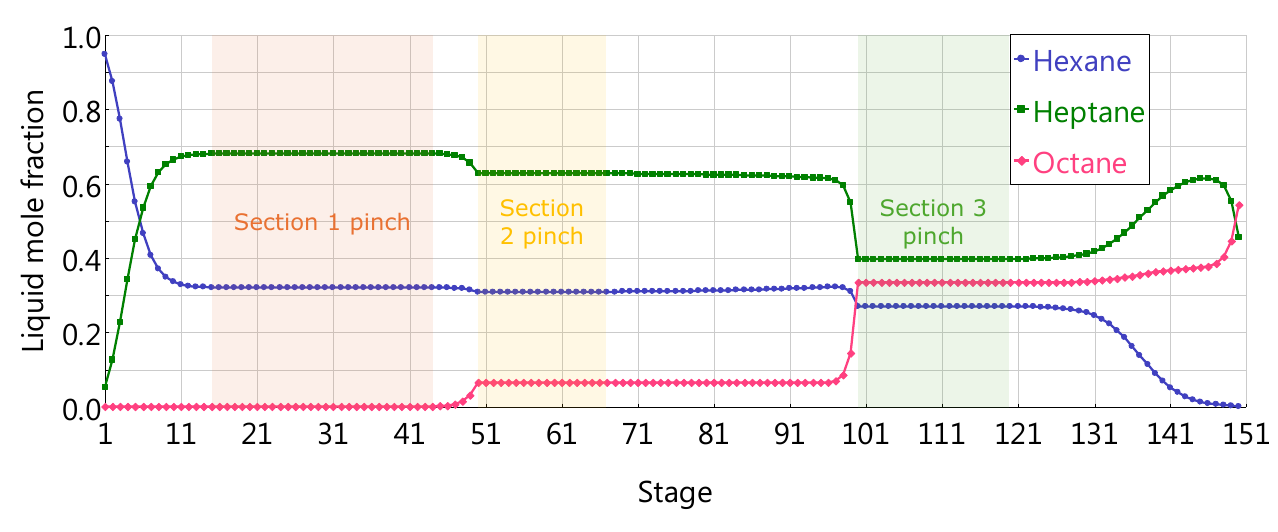}
	\caption{The liquid composition profile retrieved from Aspen Plus at the true minimum reflux ratio of $R_{\min}=2.145$. It is clear that the pinch in $\SEC_1$ is a saddle point and is located inside the column section, whereas the pinches in $\SEC_2$ and $\SEC_3$ are both unstable nodes and are located at the top of their column sections. Note that the colors of these pinch zones match with their pinch simplices drawn in Figure \ref{fig_case1scenario1}.}
	\label{fig_aspen}
\end{figure}

We validate the minimum reflux ratio obtained from our shortcut method using rigorous Aspen Plus simulation. Each column section contains 50 equilibrium stages, much larger than what are needed for this paraffin separation task. This is to ensure that the true minimum reflux condition is achieved. It turns out that the minimum reflux ratio obtained from our shortcut method is less than 1\% different compared to true minimum reflux ratio ($R_{\min}=2.145$) obtained from rigorous Aspen Plus simulation. Also, the liquid composition profile inside this two-feed column at minimum reflux, as shown in Figure \ref{fig_aspen}, exactly follows the behavior of liquid composition trajectory bundle of a pinch simplex (see Figure \ref{fig_simplex}). For more details, readers are directed to review Sections 3.4 and 4.2 of Jiang et al.\cite{jiang}. Specifically, since the distillate product is free of n-octane, the liquid composition profile $\mathbf{x}_n$ (where stage number $n$ is numbered from top to bottom) starting from the distillate product ($n=0$) must lie on the hyperplane $z_{1}^{\SEC_1}(\mathbf{x}) = 0$ until it reaches a (saddle) pinch, which corresponds to vertex $Z_2$ of the pinch simplex and lies inside $\SEC_1$. Below this pinch, the liquid composition profile continues along the hyperplane $z_{3}^{\SEC_1}(\mathbf{x}) = 0$ until it reaches the lower end of $\SEC_1$, which is connected to the top of $\SEC_2$. It turns out this is where the pinch zone lies for $\SEC_2$. Since the pinch is an unstable node when moving downward along the column, the liquid composition profile moves away from the pinch until it reaches the lower end of $\SEC_2$. Again, the pinch zone of $\SEC_3$ is located at the top of the section, from which the composition profile follows its trajectory inside the pinch simplex and heads toward the stable node ($Z_1$) until it reaches the bottoms product composition. It is worth noting that, while the n-hexane composition is small in the bottoms product (0.1 mol\%), it is not negligible. Thus, although the liquid composition profile inside $\SEC_3$ may appear to be approaching to the saddle point pinch, it never actually reaches the saddle pinch, which can be seen from Figure \ref{fig_aspen}. 

Next, using the same two-feed column example, we will examine the prevailing modeling heuristics that (1) a MFMP column can be decomposed into a series of simple columns with exactly one feed and two products, and (2) the actual minimum reflux ratio of the original MFMP column is simply the largest minimum reflux ratio value determined for all decomposed simple columns (which can be determined by applying the classic Underwood method \cite{underwood,underwood2}). According to column decomposition, the two-feed column of Figure \ref{fig_case1} is modeled as two simple columns, with one having $\Feed_1$ as the feed stream and consisting of $\SEC_1$ and $\SEC_2$, whereas the other with $\Feed_2$ as the feed stream and consisting of $\SEC_2$ and $\SEC_3$. In this case, it turns out that the largest minimum reflux ratio of the two decomposed simple columns is identified as $2.618$, which is significantly higher than the true minimum reflux ratio. In other words, the column decomposition approach overestimates the true minimum reflux in this example.

\begin{figure}[!ht]
	\centering
	\includegraphics[scale=0.7]{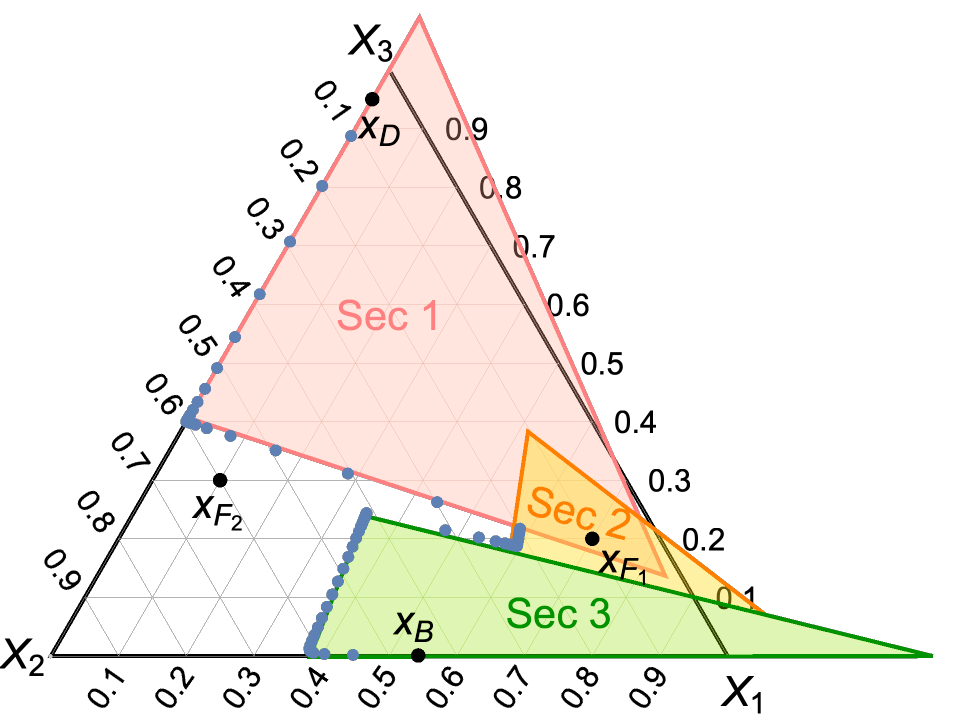}
	\caption{For the case where the upper feed is less volatile than the lower feed, the pinch simplex diagram at calculated minimum reflux ratio of $R_{\min}=1.683$. The blue dots indicate the liquid composition profile at $R = 1.738$, which is the minimum reflux ratio predicted by rigorous Aspen Plus simulation.}
	\label{fig_twofeedsimplex2}
\end{figure}

Now, we consider the second scenario where the locations of the two feed streams are switched. In other words, the upper feed $\Feed_1$ is less volatile than the lower feed $\Feed_2$. The distillate and bottoms product specifications remain unchanged. Using Algorithms \ref{algo1} through \ref{algo3}, we determine that the minimum reflux ratio of this new arrangement is $R_{\min}=1.683$, at which the lower feed $\Feed_2$ controls the separation. This can be visualized from the pinch simplex diagram of Figure \ref{fig_twofeedsimplex2}, where sections $\TOPF_2$ (i.e., $\SEC_2$) and $\BOTF_2$ (i.e., $\SEC_3$) share a common boundary, indicating that $\gamma^{\TOPF_2}_3 = \gamma^{\BOTF_2}_2 = \rho_{2,\Feed_2}$ is the binding constraint. Rigorous Aspen Plus simulation shows that the true minimum reflux ratio is $1.738$. Thus, our shortcut model gives an accurate estimation of the minimum reflux ratio with a 3\% relative difference compared to the true minimum reflux ratio. Furthermore, if we adopt the column decomposition method, we would end up with a ``minimum reflux ratio'' that is as high as 19.714, which is almost 11.3 times as large as the true minimum reflux ratio! Clearly, designing or operating the MFMP column based on incorrect minimum reflux ratio will lead to tremendous capital and operating costs.

By examining the two scenarios, we find that the optimal feed arrangement does not necessarily follow any particular pattern based on its temperature. Intuitively, one might think that, in order to reduce energy consumption (i.e., reflux ratio), feed streams should be placed according to their temperatures. Specifically, a common belief is that a high-temperature feed should be placed closer to the bottom of the column than a low-temperature feed. However, it turns out that, despite achieving the same product flow rates and purities, the minimum reflux ratio in the first scenario ($R_{\min}=2.162$) is much higher than that in the second scenario ($R_{\min}=1.683$)! This finding matches with the observation first made by Levy and Doherty \cite{levy}. And here we provide the first systematic analysis of the contradictions to the common belief that a high-temperature feed should be placed below a low-temperature feed. Practitioners should examine carefully the optimal feed arrangement when designing their columns. In this regard, our shortcut model and minimum reflux calculation method allow practitioners to obtain a quick and reliable screening of the optimal feed arrangement.

\subsection{Example 2: A One-Feed, Two-(Side)Product Column}

\begin{figure}[!ht]
	\centering
	\includegraphics[scale=1.0]{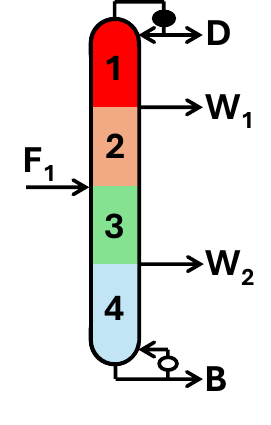}
	\vspace{-1em}
	\caption{A one-feed column with two sidedraw streams considered in Example 2. The colors of the pinch simplices match with those in Figure \ref{fig_case2}.}
	\label{fig_case2}
\end{figure}

In this example, we study a distillation column separating a ternary mixture of n-hexane, n-heptane, and n-octane with one feed stream and two sidedraw product streams, as shown in Figure \ref{fig_case2}. When there is only one feed stream and both sidedraw products are withdrawn as saturated liquids, there is a common belief in the literature (e.g., Sugie and Lu \cite{sugie}, Glinos and Malone \cite{glinos}) that $\Feed_1$ will always be ``controlling'' the separation at minimum reflux. This assumption originates from the observation of the McCabe-Thiele diagram for binary distillation with saturated liquid sidedraws, where the operating lines for sections above $\Feed_1$ continuously decreases from the top of the column to $\Feed_1$, and the operating lines for sections below $\Feed_1$ continuously increases from the bottom of the column to $\Feed_1$. 

To verify if this result can be generalized to multicomponent distillation, we present this example where $\Feed_1$ is a saturated liquid stream containing 30 mol/s of n-hexane (Component 3), 40 mol/s of n-heptane (Component 2), and 30 mol/s of n-octane (Component 1). The distillate stream contains 24 mol/s of n-hexane, 6 mol/s of n-heptane and negligible amount of n-octane, whereas the bottoms product contains 20 mol/s of n-octane and no n-hexane or n-heptane. The upper sidedraw $\Sidedraw_1$, which is located above $\Feed_1$, is a saturated liquid stream with 6 mol/s of n-hexane and 24 mol/s of n-heptane. The lower sidedraw $\Sidedraw_2$ is also a saturated liquid stream with 10 mol/s of n-heptane and 10 mol/s of n-octane. Once $\{d_i^{\SEC_k}\}_{i,k}$ are determined, we determine that $\mathcal{I}_{\Sidedraw_1} = \{2\}, \mathcal{I}_{\Feed_1} = \{2,3\}, \mathcal{I}_{\Sidedraw_2} = \{2,3\}$.

\begin{figure}[!ht]
	\centering
	\includegraphics[trim=0 0 0 11cm,clip,scale=0.7]{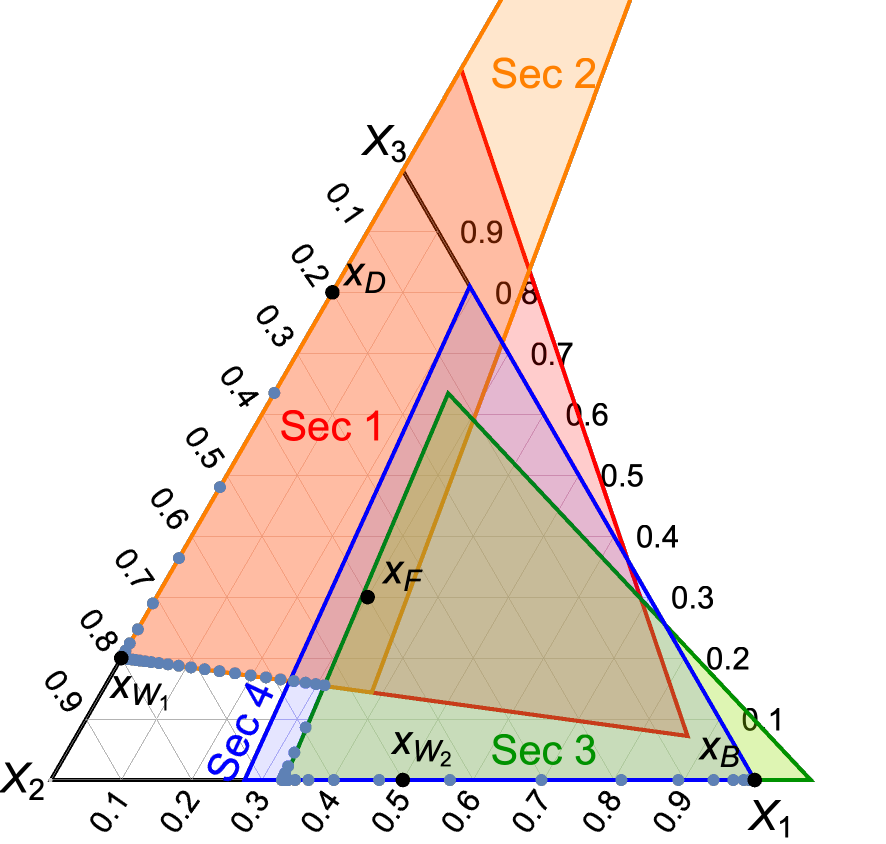}
	\caption{The pinch simplex diagram for a one-feed, two-(side)product column example at minimum reflux ($R_{\min}=2.693$), along with the liquid composition profile at the minimum reflux of $R_{\min}=2.668$ identified by Aspen Plus (also see Figure \ref{fig_aspen2}). The exact pinches compositions in $\SEC_1$ through $\SEC_4$ are $Z_2$ (associated with pinch root $\gamma_p^{\SEC_1} = \gamma_2^{\SEC_1} \in (\alpha_1, \alpha_2)$), $Z_2$ (associated with pinch root $\gamma_p^{\SEC_2} = \gamma_2^{\SEC_2} \in (\alpha_1, \alpha_2)$), $Z_2$ (associated with pinch root $\gamma_p^{\SEC_3} = \gamma_2^{\SEC_3} \in (\alpha_2, \alpha_3)$), and $Z_1$ ($\gamma_p^{\SEC_4} = \gamma_1^{\SEC_4} \in (\alpha_1, \alpha_2)$), respectively. Therefore, $\mu_2^{\SEC_1} = \mu_2^{\SEC_2} = \mu_3^{\SEC_3} = \mu_2^{\SEC_4} = 1$.}
	\label{fig_twosdsimplex}
\end{figure}

By applying Algorithms \ref{algo1} through \ref{algo3}, we determine that the minimum reflux ratio is $R_{\min}=2.693$, which is less than 1\% different compared to rigorous Aspen Plus simulation minimum reflux ratio result of $2.668$. From the minimum reflux pinch simplex diagram of Figure \ref{fig_twosdsimplex}, we can see that sidedraw $\Sidedraw_1$ actually controls the separation at minimum reflux. Specifically, the minimum reflux occurs when $z_2^{\TOPW_1}(\mathbf{x}_{\Sidedraw_1}) = z_2^{\BOTW_1}(\mathbf{x}_{\Sidedraw_1})=0$, indicating that $\gamma_3^{\TOPW_1} = \gamma_3^{\BOTW_1} = \rho_{2,\Sidedraw_1}$. This is a consequence of Equation \eqref{eqn_feasibility4}, which requires that the sidedraw composition $\mathbf{x}_{\Sidedraw_1}$ must lie on the liquid composition profile, and thus must not reside outside of the pinch simplices associated with $\SEC_1$ and $\SEC_2$. If the reflux ratio drops below this minimum threshold, the sidedraw composition $\mathbf{x}_{\Sidedraw_1}$ no longer resides within the pinch simplices, as their boundaries $z_3^{\TOPW_1}(\mathbf{x}) = 0$ and $z_3^{\BOTW_1}(\mathbf{x}) = 0$ would move toward $X_3$ (pure n-hexane).

\begin{figure}[!ht]
	\centering
	\includegraphics[width=\textwidth]{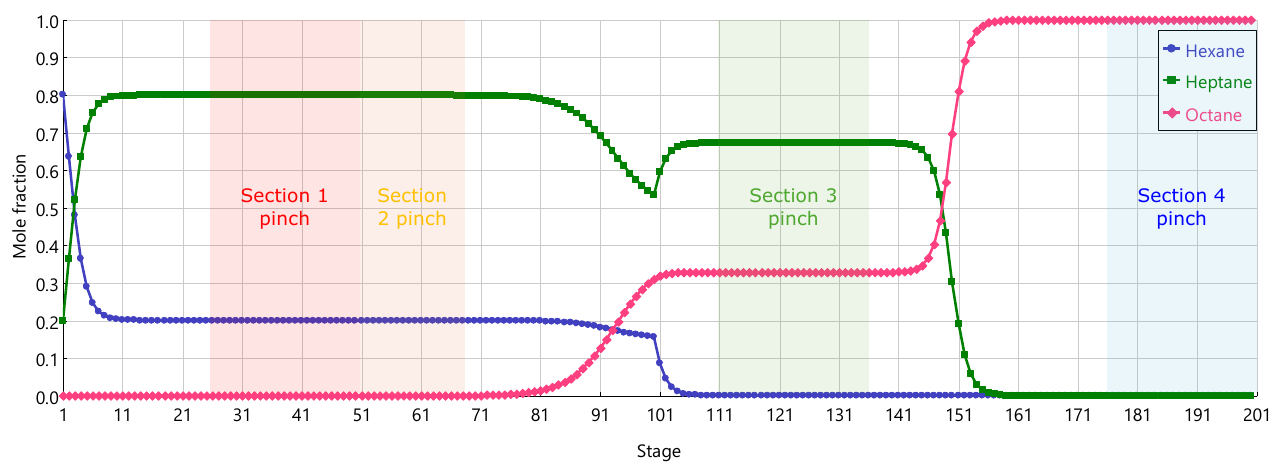}
	\vspace{-2em}
	\caption{The liquid composition profile retrieved from Aspen Plus at the true minimum reflux ratio $R_{\min}=2.668$. Each of the four column sections are given 50 equilibrium stages. The pinch zones in $\SEC_1$ through $\SEC_4$ are located respectively at the bottom, at the top, within, and at the bottom of their corresponding column sections. Note that the colors of these pinch zones match with their pinch simplices drawn in Figure \ref{fig_twosdsimplex}.}
	\label{fig_aspen2}
\end{figure}

Now, to see what happens when we insist $\Feed_1$ to control the separation at minimum reflux, we relax the feasibility constraints in Algorithms \ref{algo1} through \ref{algo3} by ignoring Equation \eqref{eqn_feasibility4}. This gives a ``minimum reflux ratio'' of $2.533$, which is lower than the true minimum reflux ratio. As a result, we have provided a counterexample to the common belief that the feed stream always controls the minimum reflux operation when sidedraws are taken as saturated liquid streams. Without incorporating the constraints related to sidedraws, one may completely ignore the possibility that a sidedraw could control the separation at minimum reflux and thus will obtain an incorrect minimum reflux ratio value that causes infeasible separation. To the best of our knowledge, this work is the first in deriving these sidedraw-related constraints and incorporating them in an algorithmic framework to calculate the true minimum reflux ratio for MFMP columns. Furthermore, we remark that our proposed minimum reflux calculation method is a generalized framework that is not limited to single-feed columns in saturated liquid sidedraws.

\subsection{Example 3: A Two-Feed, One-(Side)Product Column} \label{sec:4compcase}

\begin{figure}[!ht]
	\centering
	\includegraphics[scale=0.7]{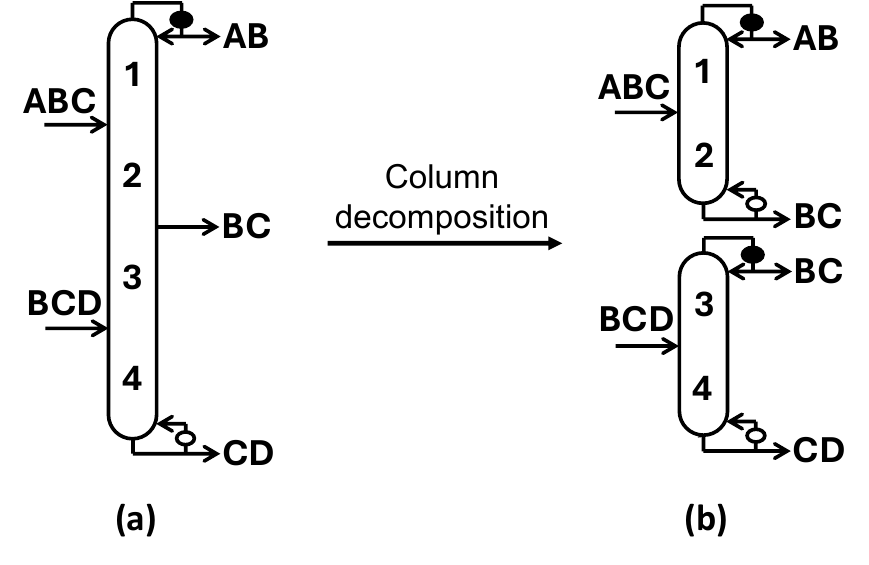}
	\caption{(a) A MFMP column for quaternary separation; (b) the decomposed version of (a).}
	\label{fig_surrogate}
\end{figure}

In the third example, we study a MFMP column drawn in Figure \ref{fig_surrogate}a that separates n-hexane (component A or 4), n-heptane (component B or 3), n-octane (component C or 2), and n-nonane (component D or 1). The relative volatilities with respect to nonane is are $(\alpha_4,\alpha_3,\alpha_2,\alpha_1) = (12.332, 5.361, 2.300, 1)$. Such a MFMP column are very common in multicomponent distillation configurations \cite{shahmatrix}, and it can be obtained by consolidating two simple columns and merging the common product stream BC into a sidedraw stream. For this MFMP column, the upper feed $\Feed_1$ (ABC) is a saturated vapor stream with 30 mol/s of n-hexane (component A), 30 mol/s of n-heptane (component B), and 40 mol/s of n-octane (component C), whereas the lower feed $\Feed_2$ (BCD) is a saturated liquid stream with 40 mol/s of n-heptane, 30 mol/s of n-octane, and 30 mol/s of n-nonane (component D). The sidedraw $\Sidedraw_1$ (BC) is in saturated liquid state. In terms of product specifications, we require that the most volatile component A must be completely recovered in the distillate stream, whereas the least volatile component D must be completely recovered in the bottoms product. The distributions of intermediate components B and C in product streams, on the other hand, are flexible. Depending on what the distributions are, component B in $\SEC_2$ and component C in $\SEC_3$ could have either net material upward or downward flow. For instance, when component B in $\SEC_2$ has net material upward flow ($d_{3}^{\SEC_2} > 0$), some of component B coming from the lower feed $\Feed_2$ will travel all the way to the top of the column and be produced as distillate. Because of this, the pinch root $\gamma_p^{\SEC_2} = \gamma_3^{\SEC_2}$ may lie in either $(\alpha_2,\alpha_3)$ or $(\alpha_3, \alpha_4)$. Similarly, the pinch root $\gamma_p^{\SEC_3} = \gamma_2^{\SEC_3}$ may lie in either $(\alpha_1,\alpha_2)$ or $(\alpha_2, \alpha_3)$. As a result, we only need one binary variable, $\mu_{3}^{\SEC_2}$, to indicate whether $\gamma_p^{\SEC_2}$ lies in $(\alpha_2,\alpha_3)$ or not. Similarly, we also only need one binary variable, $\mu_{2}^{\SEC_3}$, to indicate whether $\gamma_p^{\SEC_3}$ lies in $(\alpha_1,\alpha_2)$ or not. 

Furthermore, since $\gamma_3^{\SEC_2} \in (\alpha_2, \alpha_4)$ and $\gamma_2^{\SEC_3} \in (\alpha_1, \alpha_3)$, singularity issue might arise when implementing Equation \eqref{eqn_characteristic} in the optimization model when pinch root $\gamma_3^{\SEC_2}$ takes the value $\alpha_3$ and/or when pinch root $\gamma_{2}^{\SEC_3}$ takes the value $\alpha_2$. To avoid the singularity issue, we reformulate Equation \eqref{eqn_characteristic} by multiplying both sides of it by the bound factor (e.g., $(\alpha_3 - \gamma_3^{\SEC_2})$ for $V^{\SEC_2}$) followed by performing partial fraction decomposition. For example, the $V^{\SEC_2}$ expression can be reformulated as:
\begin{equation*}
	\begin{split}
	V^{\SEC_2}(\alpha_3 - \gamma_3^{\SEC_2}) &= (\alpha_3 - \gamma_3^{\SEC_2}) \frac{\alpha_2 d_{2}^{\SEC_2}}{\alpha_2 - \gamma_3^{\SEC_2}} + \alpha_3 d_{3}^{\SEC_2} \\
	& = \alpha_2 d_{2}^{\SEC_2} + (\alpha_3-\alpha_2) \frac{\alpha_2 d_{2}^{\SEC_2}}{\alpha_2 - \gamma_3^{\SEC_2}} + \alpha_3 d_{3}^{\SEC_2},
	\end{split}
\end{equation*}
which we note that $d_1^{\SEC_2} = d_4^{\SEC_2} =0$ because n-hexane is completely recovered in the distillate stream and n-nonane is completely recovered in the bottoms stream, as shown in Figure \ref{fig_surrogate}. Similarly, we can reformulate the $V^{\SEC_3}$ expression using this technique. With this, we can safely bound $\gamma_3^{\SEC_2} \in (\alpha_2,\alpha_4)$ and $\gamma_2^{\SEC_3} \in (\alpha_1,\alpha_3)$ witout concerning about the singularity issue. In Appendix C, we provide all the equations and constraints needed to determine the optimal distribution of intermediate components to minimize the reboiler vapor duty requirement (i.e., $V^{\SEC_4}$) for this MFMP column. The resulting optimization model, which is a mixed-integer nonlinear program (MINLP), is solved to global optimality within 15 CPU seconds in a Dell Precision 7865 workstation (equipped with 128 GB RAM and AMD Ryzen Threadripper PRO 5975WX 32-Cores 3.6 GHz processor) using global solver BARON 24.3 \cite{baron} via GAMS 46.5. In Appendix \ref{appendixC}, we provide the complete MINLP formulation implemented in GAMS. The lowest possible minimum reboiler vapor duty $V_{\text{reb},\min}$ (i.e., minimum vapor duty in $\SEC_4$) is determined to be 71.87 mol/s. And the corresponding optimal product distributions are summarized in Table \ref{table1}.

\begin{table}[ht]
\centering
 \begin{tabular}{|c|c|c|} 
	\hline
	Stream & Label in& Flow rate of component \\
	&  Figure \ref{fig_surrogate}a & A, B, C, D (mol/s) \\
	\hline
	Distillate & AB & 30, 14.23, 0, 0 \\ 
	\hline
	Sidedraw & BC  & 0, 55.77, 48.94, 0 \\
	\hline
	Bottoms & CD & 0, 0, 21.06, 30\\
	\hline
\end{tabular}
\caption{Optimal distributions of n-hexane, n-heptane, n-octane, and n-nonane in all product streams at $V_{\text{reb},\min} = 71.87$ mol/s. One can determine that $d_4^{\SEC_2}=0$ and $d_3^{\SEC_2} = -15.77 \text{ mol/s} <0$, and thus $\mu_{4}^{\SEC_2} =1$. Meanwhile, $d_2^{\SEC_3} = 8.94 \text{ mol/s} >0$ and $d_1^{\SEC_3} = 0$, and thus $\mu_{2}^{\SEC_3} = 1$.}
\label{table1}
\end{table}

We verify this result by performing exhaustive sensitivity analysis using Aspen Plus. The lowest reboiler vapor duty requirement that satisfies product requirements is found to be $77.9$ mol/s, which is within 8.3\% relative difference compared to the MINLP result. The associated n-heptane and n-octane flow rates in product streams also match very well with the results shown in Table \ref{table1}. This validates the accuracy and computationally efficiency of global optimization framework based on the shortcut model. Moreover, we remark that the global optimization algorithm does more than just finding the minimum energy requirement of a MFMP column and its corresponding product distributions. For example, there has been a lingering question among the distillation community of whether all n-heptane can be recovered from the distillate product in this MFMP column. We can easily answer questions like this by modifying the relevant variable bounds and/or by adding/deactivating related constraints in the MINLP formulation. In this case, by introducing a new constraint $d_{3}^{\SEC_1} = f_{3,\Feed_1} + f_{3,\Feed_2}$ into the MINLP formulation, the resulting optimization problem turns out to be infeasible. Thus, we conclude that it is impossible to recover all the n-heptane in the distillate product. Rigorous Aspen Plus simulation also confirms that some n-heptane is always drawn from the sidedraw no matter how much vapor is generated at the reboiler.

\begin{table}[!ht]
	\centering
	 \begin{tabular}{|c|c|c|c|} 
		\hline
		Stream & Label in & Flow rate of component \\
		 &  Figure \ref{fig_surrogate}a & A, B, C, D (mol/s) \\
		\hline
		Distillate & AB & 30, 40, 0, 0  \\ 
		\hline
		Sidedraw & BC & 0, 30, 40, 0  \\
		\hline
		Bottoms & CD & 0, 0, 30, 30 \\
		\hline
	\end{tabular}
	\caption{Component molar flow rates (arranged as n-hexane, n-heptane, n-octane, and lastly n-nonane) of all product streams.} \label{table2}
\end{table}

Lastly, using this MFMP column as an example, we illustrate why the column decomposition method shown in Figure \ref{fig_surrogate} fails to calculate the true minimum reflux ratio. The product flow rates and compositions in this example have already been specified and are listed in Table \ref{table2}. The minimum reflux ratio (which is also achieved when reboiler vapor duty is minimized as product flow rates are fixed) can be calculated using Algorithm \ref{algo1} or approach. In particular, it is worth mentioning that the resulting optimization program is solved to global optimality instantaneously during preprocessing step. Both approaches give the same minimum reflux ratio of $R_{\min}=2.002$, which is only 0.1\% different from the minimum reflux ratio of $2.000$ predicted by Aspen Plus simulation. Furthermore, this is achieved when sidedraw BC controls the minimum reflux condition. Meanwhile, the column decomposition method, which calculates the minimum reflux ratio of two simple columns as shown in Figure \ref{fig_surrogate}b using the classic Underwood method, yields a ``minimum reflux ratio'' of $1.806$, which is significantly lower than the true minimum reflux ratio. In other words, if the column operates at $R = 1.806$, the desired separation can never be achieved. 

There are two main reasons why column decomposition technique fails in this example. First, from Table \ref{table2}, one can calculate that the component distillate flow for n-heptane (B) is greater than n-heptane flow rate in the upper feed $\Feed_1$. This means that some of n-heptane in the distillate must come from the lower feed $\Feed_2$. Likewise, since the component bottoms flow for n-octane (C) is greater than the n-octane flow rate in the lower feed, some of n-octane in the bottoms must come from the upper feed $\Feed_1$. Therefore, in this MFMP column, components with intermediate relative volatilities do not follow the same flow pattern when the MFMP column is decomposed into two single columns. The second reason is that, as the original MFMP column is decomposed into two simple columns, we lose the possibility that stream BC may control the minimum reflux. Therefore, we must consider the entire MFMP column as a whole when modeling its separation performance and determining its minimum reflux condition, using methods such as the one presented in this work.

\section{Conclusion}
In this paper, we introduce the mathematical formulation that incorporates the model developed in the first article of the series \cite{jiang} to determine the minimum reflux condition of MFMP columns for multicomponent distillation. When the full product specifications are given, an algorithmic procedure is developed to automatically determine the minimum reflux ratio or minimum reboiler vapor duty requirement. When some of the product specifications are not given to users a priori, an optimization model can be developed as an MINLP to simulatenously identify the minimum reflux ratio and the corresponding optimal product distributions. We present the use of both approaches to analyze the minimum reflux behavior of MFMP columns. In all case studies, our minimum reflux ratio results matches very well with rigorous Aspen Plus simulation results.

In addition to validating the accuracy and usefulness of our proposed algorithmic and optimization frameworks, the second aim of these case studies is to reexamine some of the well-accepted design heuristics and modeling assumptions the distillation community has been relying on regarding how MFMP columns should be designed and operated. It turns out that some of these heuristics and assumptions need to be rewritten. In Example 1, we show a counterexample where placing a colder feed stream above a warmer feed stream, which follows the temperature profile within the column, actually leads to a higher minimum vapor duty requirement than if the feed stream locations are reversed. Thus, we must analyze all possible permutations of relative feed locations to determine the optimal feed stream arrangement. Our shortcut based approach is particularly suitable for analyses like this compared to rigorous process simulations which can be quite time consuming to perform, especially as the number of feed streams and/or sidedraw streams increases.

Another key finding is that decomposing a MFMP column into multiple simple columns and taking the largest individual minimum reflux ratios of each decomposed column using the classic Underwood method is not the correct approach to determine the minimum reflux ratio of the original MFMP column. In fact, such column decomposition approach can lead to minimum reflux ratio values that significantly deviate from the true minimum reflux ratio. On the other hand, our shortcut based approach considers the entire MFMP column as a whole, which is needed for accurately estimating the true minimum reflux ratio.

Finally, when a distillation column has one or more sidedraw streams, one of the sidedraw streams can control the separation at minimum reflux, even when they are all withdrawn as saturated liquid streams. This possibility has often been overlooked by the distillation community in the past due to the lack of fundamental understanding and systematic tools to model how sidedraws affect the minimum reflux operation of a multicomponent distillation column. The mathematical model and algorithms developed in this series have filled this gap, thus allowing practitioners to conduct rigorous, accurate analysis of columns with sidedraws for the first time. Overall, we believe that these new findings and insights are helpful in synthesizing and operating energy-efficient, cost-competitive, and intensified MFMP columns for multicomponent distillation.

\section*{Acknowledgment}
The information, data, or work presented herein was funded in part by the start-up fund from College of Engineering, Archiecture, and Technology at Oklahoma State University.

\section*{Data Availability and Reproducibility Statement}
The liquid composition data obtained from rigorous Aspen Plus simulations are transformed into the equilateral triangular coordinate as shown in Figures \ref{fig_case1scenario1}, \ref{fig_twofeedsimplex2}, and \ref{fig_twosdsimplex} to be visualized. The transformed liquid composition data are provided in the Supplementary Material. The procedure and specifications used to produce the pinch simplices for all column sections are provided in Algorithms \ref{algo1} through \ref{algo3} within the article. The optimal product distribution results in Table \ref{table1} are obtained by the MINLP formulation provided in Appendix \ref{appendixC}. 

\bibliographystyle{aichej}

\bibliography{references}

\begin{thebibliography}{10}
\providecommand{\url}[1]{\texttt{#1}}
\providecommand{\urlprefix}{URL }

\bibitem{distillationbook}
Górak A, Olujić Z.
\newblock \emph{Distillation: Fundamentals and Principles}.
\newblock Elsevier Inc. 2014.

\bibitem{amo}
{US DOE Industrial Efficiency \& Decarbonization Office}.
\newblock {Manufacturing Energy and Carbon Footprints (2018 MECS)}. 2021.

\bibitem{roadmap}
{US DOE}.
\newblock Industrial Decarbonization Roadmap. 2022.

\bibitem{shahmatrix}
Shah VH, Agrawal R.
\newblock A matrix method for multicomponent distillation sequences.
\newblock \emph{AIChE Journal}. 2010;\hspace{0pt}56(7):1759--1775.

\bibitem{gmaneed}
Nallasivam U, Shah VH, Shenvi AA, Tawarmalani M, Agrawal R.
\newblock Global optimization of multicomponent distillation configurations: 1. Need for a reliable global optimization algorithm.
\newblock \emph{AIChE Journal}. 2013;\hspace{0pt}59(3):971--981.

\bibitem{nallasivam}
Nallasivam U, Shah VH, Shenvi AA, Huff J, Tawarmalani M, Agrawal R.
\newblock Global optimization of multicomponent distillation configurations: 2. Enumeration based global minimization algorithm.
\newblock \emph{AIChE Journal}. 2016;\hspace{0pt}62(6):2071--2086.

\bibitem{pireview}
Jiang Z, Agrawal R.
\newblock Process intensification in multicomponent distillation: A review of recent advancements.
\newblock \emph{Chemical Engineering Research and Design}. 2019;\hspace{0pt}147:122--145.

\bibitem{pireview2}
Jiang G Zheyuand Madenoor~Ramapriya, Tawarmalani M, Agrawal R.
\newblock Process intensification in multicomponent distillation.
\newblock \emph{Chemical Engineering Transactions}. 2018;\hspace{0pt}69:841--846.

\bibitem{misconception}
Agrawal R, Tumbalam~Gooty R.
\newblock Misconceptions about efficiency and maturity of distillation.
\newblock \emph{AIChE Journal}. 2020;\hspace{0pt}66(8):e16294.

\bibitem{heatpump}
Agrawal R, Yee TF.
\newblock Heat Pumps for Thermally Linked Distillation Columns: An Exercise for Argon Production from Air.
\newblock \emph{Industrial \& Engineering Chemistry Research}. 1994;\hspace{0pt}33(11):2717--2730.

\bibitem{akash}
Nogaja A, Tawarmalani M, Agrawal R.
\newblock Distillation Electrification Through Optimal Use of Heat Pumps.
\newblock In: \emph{34th European Symposium on Computer Aided Process Engineering / 15th International Symposium on Process Systems Engineering}, edited by Manenti F, Reklaitis GV, vol.~53 of \emph{Computer Aided Chemical Engineering}, pp. 1297--1302. Elsevier. 2024;\hspace{0pt}.

\bibitem{wankat}
Wankat PC.
\newblock Multieffect distillation processes.
\newblock \emph{Industrial \& Engineering Chemistry Research}. 1993;\hspace{0pt}32(5):894--905.

\bibitem{intermediatereboiler}
Agrawal R, Herron DM.
\newblock Intermediate reboiler and condenser arrangement for binary distillation columns.
\newblock \emph{AIChE Journal}. 1998;\hspace{0pt}44(6):1316--1324.

\bibitem{prefractionation}
Agrawal R, Fidkowski ZT, Xu J.
\newblock Prefractionation to reduce energy consumption in distillation without changing utility temperatures.
\newblock \emph{AIChE Journal}. 1996;\hspace{0pt}42(8):2118--2127.

\bibitem{feedconditioning}
Agrawal R, Herron DM.
\newblock Optimal thermodynamic feed conditions for distillation of ideal binary mixtures.
\newblock \emph{AIChE Journal}. 1997;\hspace{0pt}43(11):2984--2996.

\bibitem{hmp}
Jiang Z, Madenoor~Ramapriya G, Tawarmalani M, Agrawal R.
\newblock Minimum energy of multicomponent distillation systems using minimum additional heat and mass integration sections.
\newblock \emph{AIChE Journal}. 2018;\hspace{0pt}64(9):3410--3418.

\bibitem{krishna}
{Tumbalam Gooty} R, {Chavez Velasco} JA, Agrawal R.
\newblock Methods to assess numerous distillation schemes for binary mixtures.
\newblock \emph{Chemical Engineering Research and Design}. 2021;\hspace{0pt}172:1--20.

\bibitem{jose}
{Chavez Velasco} JA, Tawarmalani M, Agrawal R.
\newblock Systematic Analysis Reveals Thermal Separations Are Not Necessarily Most Energy Intensive.
\newblock \emph{Joule}. 2021;\hspace{0pt}5(2):330--343.

\bibitem{gilliland2}
Gilliland ER.
\newblock MULTICOMPONENT RECTIFICATION.
\newblock \emph{Industrial \& Engineering Chemistry}. 1940;\hspace{0pt}32(8):1101--1106.

\bibitem{doherty}
Doherty MF.
\newblock \emph{Conceptual design of distillation systems}.
\newblock McGraw-Hill chemical engineering series. Boston: McGraw-Hill. 2001.

\bibitem{jiang}
Jiang Z, Tawarmalani M, Agrawal R.
\newblock Minimum reflux calculation for multicomponent distillation in multi-feed, multi-product columns: Mathematical model.
\newblock \emph{AIChE Journal}. 2022;\hspace{0pt}68(12):e17929.

\bibitem{cht}
Mathew TJ, Tawarmalani M, Agrawal R.
\newblock Relaxing the constant molar overflow assumption in distillation optimization.
\newblock \emph{AIChE Journal}. 2023;\hspace{0pt}69(9):e18125.

\bibitem{nlpcost}
Jiang Z, Mathew TJ, Zhang H, Huff J, Nallasivam U, Tawarmalani M, Agrawal R.
\newblock Global optimization of multicomponent distillation configurations: Global minimization of total cost for multicomponent mixture separations.
\newblock \emph{Computers \& Chemical Engineering}. 2019;\hspace{0pt}126:249--262.

\bibitem{nlpexergy}
Jiang Z, Chen Z, Huff J, Shenvi AA, Tawarmalani M, Agrawal R.
\newblock Global minimization of total exergy loss of multicomponent distillation configurations.
\newblock \emph{AIChE Journal}. 2019;\hspace{0pt}65(11):e16737.

\bibitem{baron}
Tawarmalani M, Sahinidis NV.
\newblock {A polyhedral branch-and-cut approach to global optimization}.
\newblock \emph{Mathematical Programming}. 2005;\hspace{0pt}103:225--249.

\bibitem{gauthamtwofeed}
Madenoor RG, Mohit T, Rakesh A.
\newblock Modified basic distillation configurations with intermediate sections for energy savings.
\newblock \emph{AIChE Journal}. 2014;\hspace{0pt}60(3):1091--1097.

\bibitem{shah}
Shah VH.
\newblock Energy savings in distillation via identification of useful configurations.
\newblock Ph.D. thesis, Purdue University. 2010.

\bibitem{underwood}
Underwood AJV.
\newblock Fractional Distillation of Multicomponent Mixtures.
\newblock \emph{Industrial \& Engineering Chemistry}. 1949;\hspace{0pt}41(12):2844--2847.

\bibitem{underwood2}
Underwood AJV.
\newblock Fractional distillation of multicomponent mixtures.
\newblock \emph{Chemical Engineering Progress}. 1948;\hspace{0pt}44:603--614.

\bibitem{levy}
Levy SG, Doherty MF.
\newblock Design and synthesis of homogeneous azeotropic distillations. 4. Minimum reflux calculations for multiple-feed columns.
\newblock \emph{Industrial \& Engineering Chemistry Fundamentals}. 1986;\hspace{0pt}25(2):269--279.

\bibitem{sugie}
Sugie H, Lu B.
\newblock On the determination of minimum reflux ratio for a multicomponent distillation column with any number of side-cut streams.
\newblock \emph{Chemical Engineering Science}. 1970;\hspace{0pt}25(12):1837--1846.

\bibitem{glinos}
Glinos KN, Malone MF.
\newblock Design of sidestream distillation columns.
\newblock \emph{Industrial \& Engineering Chemistry Process Design and Development}. 1985;\hspace{0pt}24(3):822--828.

\end{thebibliography}

\appendix
\section{Parameters and Variables} \label{appendixA}
\begin{table}[H]
	\centering
	\begin{tabular}{l@{\hskip 0.1in}l@{\hskip 0.1in}}
	    \toprule
		$c$				 					& Total number of components present in the distillation column\\
		$\alpha_i$		 					& Relative volatility of component $i$ with respect to the heaviest component\\
		$d_i^{\SEC_k}$           			& Component $i$'s net material upward flow in column section $k$ \\
		$V^{\SEC_k}$             			& Total vapor flow in column section $k$ \\
		$\gamma_i^{\SEC_k}$      			& The $i^{\textrm{th}}$ root of Equation \eqref{eqn_characteristic} \\
		$p^{\SEC_k}$						& The pinch index of section $k$ \\
		$\gamma_p^{\SEC_k}$      			& The pinch root of section $k$ \\
		$\mathbf{x}_n$ 						& The vector of liquid composition leaving stage $n$\\
		$X_i$								& The composition of pure component $i$ \\
		$\mathbf{x}_{\Feed_j}$	 			& Liquid composition (or hypothetical liquid composition) composition of $j^{\text{th}}$ feed stream \\
		$\mathbf{x}_{\Sidedraw_j}$ 			& Liquid composition (or hypothetical liquid composition) of $j^{\text{th}}$ sidedraw stream \\
		$\rho_{i,\Feed_j},\, \rho_{i,\Sidedraw_j}$	& The $i^{\textrm{th}}$ root defined in Equations \eqref{eqn_feed} and \eqref{eqn_sidedraw}, respectively\\
		$V_{\Feed_j},\, V_{\Sidedraw_j}$	& Total vapor flow rate in the $j^{\text{th}}$ feed and sidedraw stream, respectively \\
		$f_{i,\Feed_j},\, f_{i,\Sidedraw_j}$& Component $i$'s flow rate in the $j^{\text{th}}$ feed and sidedraw stream, respectively \\
		$l_{m,\Feed_j},\, l_{m,\Sidedraw_j}$ & Component $m$'s flow rate in the liquid portion of the $j^{\text{th}}$ feed and sidedraw, respectively\\
		$v_{m,\Feed_j},\, v_{m,\Sidedraw_j}$ & Component $m$'s flow rate in the vapor portion of the $j^{\text{th}}$ feed and sidedraw, respectively\\
		$N_{\Feed}$							& Number of feed streams in the column\\
		$N_{\Sidedraw}$						& Number of sidedraw streams in the column\\
		$N_{\SEC}$							& Number of column sections in the column, which is equal to $N_{\Feed} + N_{\Sidedraw} + 1$\\
		$\mu_{i}^{\SEC_k}$					& Binary variable that equals 1 when $\gamma_p^{\SEC_k} \in (\alpha_{i-1},\alpha_i)$, and 0 otherwise\\
		$K_i^{\SEC_k}$						& Binary variable defined in Equation \eqref{eqn_K} and used in Equations \eqref{eqn_feasibility_K_feed}, \eqref{eqn_feasibility_K_sidedraw} and \eqref{eqn_feasibility4}\\
		\bottomrule
	\end{tabular}
\end{table}

\section{Sets and Notations} \label{appendixB}
\begin{table}[H]
	\centering
	\begin{tabular}{l@{\hskip 0.1in}l@{\hskip 0.1in}}
		\toprule
		$\mathcal{C}$     			&	$\{1,\cdots,c\}$\\
		$\Feed_j, \, \Sidedraw_j$   &	The $j^{\text{th}}$ feed and sidedraw stream (both counting from the top), respectively\\
		$\mathcal{I}_{\Feed_j}$		&	Index set defined in Equation \eqref{eqn_Ifeed} for feasibility constraints associated with the $j^{\text{th}}$ feed stream\\
		$\mathcal{I}_{\Sidedraw_j}$	&	Index set defined in Equation \eqref{eqn_ISidedraw} for feasibility constraints associated with the $j^{\text{th}}$ sidedraw\\
		$\TOPF_j$	    			&	Column section above the $j^{\text{th}}$ feed stream\\
		$\BOTF_j$	    			&	Column section below the $j^{\text{th}}$ feed stream\\
		$\TOPW_j$	    			&	Column section above the $j^{\text{th}}$ sidedraw stream, respectively\\
		$\BOTW_j$	    			&	Column section below the $j^{\text{th}}$ sidedraw stream, respectively\\
		\bottomrule
	\end{tabular}
\end{table}

\section{Optimization Model}\label{appendixC}

Here, we provide the complete formulation for identifying the minimum reboiler vapor duty and the corresponding product distributions for the column shown in Figure \ref{fig_case1}.

\textbf{Objective function:}
\begin{equation*}
\textrm{minimize } V^{\SEC_4}
\end{equation*}

\textbf{Constraints and bounds:}

1. Mass balance equations and feed and product specifications:
\begin{equation*}
\begin{split}
	& f_{i, \Feed_1} + f_{i, \Feed_2} = d_i^{\SEC_1} - f_{i, \Sidedraw_1} - d_{i}^{\SEC_4} \qquad \forall i=1,\dots,4 \\
	& d_{i}^{\SEC_2} = d_{i}^{\SEC_1} - f_{i,\Feed_1};\; d_{i}^{\SEC_3} = d_{i}^{\SEC_2} - f_{i,\Sidedraw_1}; d_{i}^{\SEC_4} = d_{i}^{\SEC_3} - f_{i,\Feed_2} \qquad \forall i=1,\dots,4 \\
	& f_{4,\Feed_1} = d_{4}^{\SEC_1} = 30;\; f_{1,\Feed_2} = -d_{1}^{\SEC_4} = 30 \\
	& d_{1}^{\SEC_1} = d_{2}^{\SEC_1} = 0;\; d_{1}^{\SEC_2} = d_{4}^{\SEC_2} = 0;\; d_{1}^{\SEC_3} = d_{4}^{\SEC_3} = 0;\; \; d_{4}^{\SEC_4} = d_{3}^{\SEC_4} = 0;\\
	& f_{3,\Feed_1} = 30;\; f_{2,\Feed_1} = 40;\; f_{1,\Feed_1} = 0;\; f_{4,\Feed_2} = 0;\; f_{3,\Feed_2} = 40;\; f_{2,\Feed_2} = 30;\;\\
	& d_{i}^{\SEC_1} \geq 0; \; f_{i,\Sidedraw_1} \leq 0;\; d_{i}^{\SEC_4} \leq 0 \qquad \forall i=1,\dots,4
	\end{split}
\end{equation*}

2. Vapor duty calculations based on Equation \eqref{eqn_characteristic}:
\begin{equation*}
\begin{split}
	V^{\SEC_1} &= \sum^4_{i=1}\frac{\alpha_i d_{i}^{\SEC_1}}{\alpha_i-\gamma_{4}^{\SEC_1}} \\
	V^{\SEC_1} &= \sum^4_{i=1}\frac{\alpha_i d_{i}^{\SEC_1}}{\alpha_i-\gamma_{3}^{\SEC_1}} \\
	V^{\SEC_2}(\alpha_3-\gamma_{3}^{\SEC_2}) &= \alpha_2 d_{2}^{\SEC_2} + (\alpha_3-\alpha_2) \frac{\alpha_2 d_{2}^{\SEC_2}}{\alpha_2 - \gamma_3^{\SEC_2}} + \alpha_3 d_{3}^{\SEC_2}\\
	V^{\SEC_2} &= \sum^4_{i=1}\frac{\alpha_id_{i}^{\SEC_2}}{\alpha_i-\gamma_{2}^{\SEC_2}} \\
	V^{\SEC_3} (\alpha_2-\gamma_{2}^{\SEC_3}) &= \alpha_2d_{2}^{\SEC_3}-(\alpha_3-\alpha_2)\frac{\alpha_3 d_{3}^{\SEC_3}}{\alpha_3-\gamma_{2}^{\SEC_3}} + \alpha_3 d_{3}^{\SEC_3}\\
	V^{\SEC_3} &= \sum^4_{i=1}\frac{\alpha_id_{i}^{\SEC_3}}{\alpha_i-\gamma_{3}^{\SEC_3}} \\
	V^{\SEC_4} &= \sum^4_{i=1}\frac{\alpha_i d_{i}^{\SEC_4}}{\alpha_i-\gamma_{1}^{\SEC_4}} \\
	V^{\SEC_4} &= \sum^4_{i=1}\frac{\alpha_i d_{i}^{\SEC_4}}{\alpha_i-\gamma_{2}^{\SEC_4}} \\
\end{split}
\end{equation*}

Note that here, only a subset of $\gamma_i^{\SEC_k}$ roots for each section need to be defined as decision variables, as others are fixed to the respective relative volatility values due to Equation \eqref{eqn_roots}. For example, $\gamma_2^{\SEC_1} = \alpha_2$ and $\gamma_1^{\SEC_1} = \alpha_1$, since the distilalte product contains no n-octane and n-nonane components ($d_2^{\SEC_1} = d_1^{\SEC_1} = 0$).

3. Defining equations for the feed and sidedraw streams:
\begin{equation*}
\begin{split}
	& V_{\Feed_1} = 100 = \sum^4_{i=1} \frac{\alpha_i f_{i,\Feed_1}}{\alpha_i-\rho_{j,\Feed_1}} \qquad \forall j=2,3\\
	& V_{\Feed_2} = 0 = \sum^4_{i=1} \frac{\alpha_i f_{i,\Feed_2}}{\alpha_i-\rho_{j,\Feed_2}} \qquad \forall j=1,2\\
	& V_{\Sidedraw_1} = 0 = \sum^4_{i=1} \frac{\alpha_i f_{i,\Sidedraw_1}}{\alpha_i-\rho_{j,\Sidedraw_1}} \qquad \forall j = 2\\
\end{split}
\end{equation*}

4. Variable bounds:
\begin{equation*}
\begin{split}
	& \gamma_{3}^{\SEC_1} \in (\alpha_2,\alpha_3);\; \gamma_{4}^{\SEC_1} \in(\alpha_3,\alpha_4) \\
	& \gamma_{2}^{\SEC_2} \in(\alpha_2,\alpha_3);\; \gamma_{3}^{\SEC_2} \in(\alpha_2,\alpha_4) \\
	& \gamma_{2}^{\SEC_3} \in (\alpha_1,\alpha_3);\; \gamma_{3}^{\SEC_3} \in(\alpha_2,\alpha_3) \\
	& \gamma_{1}^{\SEC_4} \in (\alpha_1,\alpha_2);\; \gamma_{2}^{\SEC_4} \in (\alpha_2,\alpha_3) \\
	& \rho_{1,\Feed_2} \in (\alpha_1,\alpha_2) ;\; \rho_{2,\Feed_1},\, \rho_{2,\Sidedraw_1},\, \rho_{2,\Feed_2} \in (\alpha_2,\alpha_3) \\
	& \rho_{3,\Feed_1} \in (\alpha_3,\alpha_4);\;
\end{split}
\end{equation*}

In actual GAMS implementation, we specify an absolute tolerance of $10^{-4}$ to bound the variables above (resp. below) from these lower (resp. upper) bounds. 

5. Binary variables $\mu_i$ and their bounds:
\begin{equation*}
\begin{split}
	&\mu_{3}^{\SEC_2},\, \mu_{2}^{\SEC_3} \in\{0,1\} \\
	& \alpha_2 \mu_{3}^{\SEC_2} + \alpha_3 (1-\mu_{3}^{\SEC_2}) \leq \gamma_{3}^{\SEC_2} \leq \alpha_3 \mu_{3}^{\SEC_2} + \alpha_4 (1-\mu_{3}^{\SEC_2}) \\ 
	& \alpha_1 \mu_{2}^{\SEC_3} + \alpha_2 (1-\mu_{2}^{\SEC_3}) \leq \gamma_{2}^{\SEC_3} \leq \alpha_2 \mu_{2}^{\SEC_3} + \alpha_3 (1-\mu_{2}^{\SEC_3})
\end{split}
\end{equation*}

6. Feasibility constraints based on Equations \eqref{eqn_feasibility_K_feed} and \eqref{eqn_feasibility_K_sidedraw}:
\begin{equation*}
\begin{split}
	& \gamma_{3}^{\SEC_1} \geq \rho_{2,\Feed_1} \geq \gamma_{2}^{\SEC_2};\; (1-\mu_{3}^{\SEC_2})(\gamma_{4}^{\SEC_1}-\rho_{3,\Feed_1}) \geq 0;\; (1-\mu_{3}^{\SEC_2})(\rho_{3,\Feed_1} - \gamma_{3}^{\SEC_2}) \geq 0 \\
	& \mu_{2}^{\SEC_3} (\gamma_{2}^{\SEC_3}-\gamma_{1}^{\SEC_2}) \geq 0;\; \gamma_{3}^{\SEC_3} \geq \rho_{2,\Sidedraw_1} \geq \gamma_{2}^{\SEC_2};\;  (1-\mu_{3}^{\SEC_2})(\gamma_{4}^{\SEC_3}-\gamma_{3}^{\SEC_2}) \geq 0 \\
	& \mu_{2}^{\SEC_3}(\gamma_{2}^{\SEC_3} - \rho_{1,\Feed_2}) \geq 0;\; \mu_{2}^{\SEC_3}(\rho_{1,\Feed_2} - \gamma_{1}^{\SEC_4}) \geq 0;\; \gamma_{3}^{\SEC_3} \geq \rho_{2,\Feed_2} \geq \gamma_{2}^{\SEC_4} \\
\end{split}
\end{equation*}

Note that here, we only need to include constraints related to $\rho_{2,\Sidedraw_1}$, as the sidedraw stream contains only n-heptane and n-octane. Therefore, when implementing Equation \eqref{eqn_feasibility_K_sidedraw} that involves $\rho_{1,\Sidedraw_1}$ (resp. $\rho_{3,\Sidedraw_1}$), i.e., $\mu_{2}^{\SEC_3} (\gamma_{2}^{\SEC_3}-\rho_{1,\Sidedraw_1}) \geq 0;\; \mu_{2}^{\SEC_3} (\rho_{1,\Sidedraw_1}-\gamma_{1}^{\SEC_2}) \geq 0$ (resp. $(1-\mu_{3}^{\SEC_2})(\gamma_{4}^{\SEC_3}-\rho_{3,\Sidedraw_1}) \geq 0;\; (1-\mu_{3}^{\SEC_2})(\rho_{3,\Sidedraw_1} - \gamma_{3}^{\SEC_2}) \geq 0$), we add them to eliminate $\rho_{1,\Sidedraw_1}$ (resp. $\rho_{3,\Sidedraw_1}$).

6. Additional constraints for sidedraw $\Sidedraw_1$ based on Equation \eqref{eqn_feasibility4}:
\begin{equation*}
	\begin{split}
		& \gamma_2^{\SEC_2} \leq \rho_{2,\Sidedraw_1} ;\; \gamma_{3}^{\SEC_3} \geq \rho_{2,\Sidedraw_1} \\
		& \mu_{3}^{\SEC_2}(\gamma_{3}^{\SEC_2} - \rho_{2,\Sidedraw_1}) \geq 0;\; (1-\mu_{2}^{\SEC_3})(\gamma_{2}^{\SEC_3}-\rho_{2,\Sidedraw_1}) \leq 0 \\
	\end{split}
\end{equation*}

Here, we also only need to include constraints related to $\rho_{2,\Sidedraw_1}$.

Figure \ref{fig_appendix} below summarizes the relationships between $\mu_i$ and $K_i$ variables for all column sections, which are used to write down the feasibility constraints and additional constraints based on Equations \eqref{eqn_feasibility_K_feed}, \eqref{eqn_feasibility_K_sidedraw} and \eqref{eqn_feasibility4} above.

\begin{figure}[!ht]
	\centering
	\includegraphics[width=0.8\textwidth]{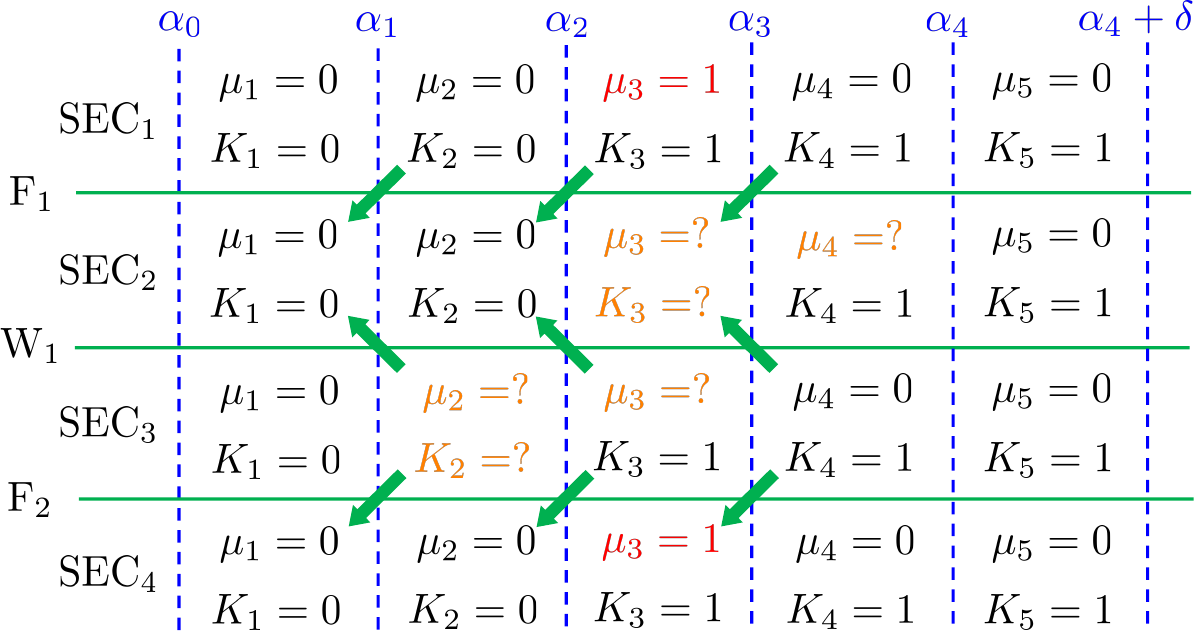}
	\caption{The relationships between $\mu_i$ and $K_i$ variables for all column sections in the two-feed, one-(side)product column example discussed in Section \ref{sec:4compcase} for quaternary separation. When the pinch root location for a specific column section is uncertain (labeled as ``?'' in the figure), it then needs to be defined as a binary variable $\mu_i$. As shown in the constraints above, in actual implementation, we replace $\mu_4^{\SEC_2}$ and $\mu_3^{\SEC_3}$ with $1-\mu_3^{\SEC_2}$ and $1- \mu_2^{\SEC_3}$, respectively. The green arrows show how the (binary) coefficients present in the feasibility constraints of Equations \eqref{eqn_feasibility_K_feed} and \eqref{eqn_feasibility_K_sidedraw} are constructed.}
	\label{fig_appendix}
\end{figure}

\end{document}